\documentclass[12pt]{amsart}

\usepackage[T1]{fontenc}
\usepackage[utf8]{inputenc}
\usepackage[english]{babel}
\usepackage{amsmath,amssymb,amscd,amsthm}
\usepackage{latexsym}
\usepackage{graphicx}
\usepackage{pgf}
\usepackage{colortbl}
\usepackage{url}
\usepackage[a4paper,twoside,left=3cm,right=2.8cm,top=3.1cm,bottom=2.3cm]{geometry}

\newtheorem{theorem}{Theorem}[section]
\newtheorem{thm}[theorem]{Theorem}

\newtheorem{corollary}[theorem]{Corollary}
\newtheorem{cor}[theorem]{Corollary}

\newtheorem{proposition}[theorem]{Proposition}
\newtheorem{lemma}[theorem]{Lemma}

\newtheorem{remark}[theorem]{Remark}

\def\e{\varepsilon}
\def\irr#1{{\rm Irr}(#1)}
\def\irrr#1#2 {\irr {#1 \mid #2}}

\newcommand{\RR}{\mathbb R}

\newcommand{\dd}{\mathrm{d}}
\newcommand{\II}{\mathrm{II}}
\newcommand{\HH}{\mathcal{H}}
\newcommand{\id}{\mathrm{Id}}
\newcommand{\oo}{\mathbf{1}}
\DeclareMathOperator{\vol}{Vol}
\DeclareMathOperator{\ent}{Ent}
\let\div\relax
\DeclareMathOperator{\div}{div}
\newcommand{\nn}[1]{\left\lVert#1\right\rVert}
\DeclareMathOperator{\var}{Var}

\def\XXint#1#2#3{{\setbox0=\hbox{$#1{#2#3}{\int}$}
     \vcenter{\hbox{$#2#3$}}\kern-.5\wd0}}

\begin{document}

\title[Remarks on the Brascamp--Lieb inequality]{Remarks on the Brascamp--Lieb inequality}
\author[D. Cordero-Erausquin, G. V. Livshyts, L. Rotem]{Dario Cordero-Erausquin, Galyna V. Livshyts, Liran Rotem}
\address{Sorbonne Universit\'e}
\email{dario.cordero@imj-prg.fr}
\address{School of Mathematics, Georgia Institute of Technology, Atlanta, GA}
\email{glivshyts6@math.gatech.edu}
\address{Technion, Israel}
\email{lrotem@technion.ac.il}

\begin{abstract}
We prove several variants of the Brascamp--Lieb inequality which are adapted
to additional symmetries, weights, and non-Euclidean radial structures. First,
an improved Brascamp--Lieb inequality for even functions yields a strengthened
Pr\'ekopa--Leindler theorem for even strongly log-concave functions. We also
record a Legendre-transform proof of a generalized log-Sobolev inequality for
log-concave functions. Further results include a Brascamp--Lieb inequality
with a concave weight, an unconditional second-order Poincar\'e inequality, a
rotation-trick improvement for arbitrary even log-concave measures, and a
Poincar\'e-type inequality associated with smooth norms.
\end{abstract}
\maketitle
\setcounter{tocdepth}{1}
\tableofcontents

\section{Introduction}

The Brascamp--Lieb inequality \cite{BrascampLiebExtensions1976}
is a fundamental functional inequality with deep connections to convex
geometry. Let $V:\RR^{n}\to\RR$ be a $C^{2}$ function such that
$\nabla^{2}V\succ0$ at every point. Let $\mu$ be the probability
measure on $\RR^{n}$ with density $Z^{-1}e^{-V}$, where
$Z=\int e^{-V}$ (unless explicitly stated otherwise, all integrations
are over $\RR^{n}$). For a $C^{1}$ function $f:\RR^{n}\to\RR$ with
$\int |f|\dd\mu<\infty$, define
\[
\var_{\mu}f=\int f^{2}\dd\mu-\left(\int f\dd\mu\right)^{2}.
\]
The Brascamp--Lieb inequality then states that for such a function $f$ we have
\begin{equation}
\var_{\mu}f\le
\int\left\langle
\left(\nabla^{2}V\right)^{-1}\nabla f,\nabla f
\right\rangle\dd\mu .
\label{eq:BL}
\end{equation}
In particular, if $\nabla^{2}V\succeq\alpha\id$ for some $\alpha>0$,
then
\[
\var_{\mu}f\le\frac{1}{\alpha}\int |\nabla f|^{2}\dd\mu .
\]
For example, the standard Gaussian measure $\gamma$, with density
$\frac{\dd\gamma}{\dd x}=Z^{-1}e^{-|x|^{2}/2}$, satisfies a Poincar\'e
inequality with constant $1$; this was known as a folklore result decades before the
work of Brascamp and Lieb. In this
note we prove new variants of several functional inequalities related to
the Brascamp--Lieb inequality.

We begin in Section \ref{sec:PL}, which concerns the Pr\'ekopa--Leindler
inequality \cite{PrekopaLogarithmic1971,PrekopaLogarithmic1973,LeindlerCertain1972}.
This inequality states that for all $0\le t\le1$ and for all measurable
functions $f,g,h:\RR^{n}\to(-\infty,\infty]$ satisfying
\[
h((1-t)x+ty)\le(1-t)f(x)+tg(y)
\]
for all $x,y\in\RR^{n}$, one has
\begin{equation}
\int e^{-h}\ge
\left(\int e^{-f}\right)^{1-t}
\left(\int e^{-g}\right)^t .
\label{eq:PL}
\end{equation}
The Pr\'ekopa--Leindler inequality is a functional version of the
Brunn--Minkowski inequality: For Borel sets $A,B\subseteq\RR^{n}$, take
$e^{-f}=\oo_A$, $e^{-g}=\oo_B$, and
$e^{-h}=\oo_{(1-t)A+tB}$ to obtain
\[
\vol((1-t)A+tB)\ge \vol(A)^{1-t}\vol(B)^t .
\]

The original paper \cite{BrascampLiebExtensions1976} of Brascamp and Lieb
studies both extensions of the Pr\'ekopa--Leindler inequality and the
Brascamp--Lieb inequality, but these two directions are treated rather
independently. Later, Bobkov and Ledoux \cite{BobkovLedouxBrunnMinkowski2000}
proved the Brascamp--Lieb inequality from the Pr\'ekopa--Leindler
inequality. To sketch the connection, fix convex functions
$f_0,f_1:\RR^n\to(-\infty,\infty]$ and consider
\[
f_t(z)=\inf_{z=(1-t)x+ty}
\left((1-t)f_0(x)+t f_1(y)\right).
\]
The Pr\'ekopa--Leindler inequality implies that
\begin{equation} \label{eq:alpha-integral}
\alpha(t)=\int e^{-f_t}
\end{equation}
is log-concave, i.e. $\log \alpha$ is a concave function. In particular this implies that 
$(\log \alpha)^{\prime\prime}(t) \le 0$, which yields the Brascamp--Lieb inequality.

In Section \ref{sec:PL} we reverse this line of reasoning. We first prove
a strengthened Brascamp--Lieb inequality for even functions, and then
deduce the following strengthened Pr\'ekopa--Leindler inequality.

\begin{thm}\label{thm:even-PL}
Let $f,g:\RR^{n}\to(-\infty,\infty]$ be even functions such that
$f-c|x|^{2}$ and $g-c|x|^{2}$ are convex for some $c>0$. Fix
$0\le t\le1$, and assume that a measurable function
$h:\RR^{n}\to\RR$ satisfies
\[
h((1-t)x+ty)\le(1-t)f(x)+tg(y)
\]
for all $x,y\in\RR^{n}$. Then
\[
\int e^{-h}e^{\frac{c}{2}|x|^{2}}
\ge
\left(\int e^{-f}e^{\frac{c}{2}|x|^{2}}\right)^{1-t}
\left(\int e^{-g}e^{\frac{c}{2}|x|^{2}}\right)^t .
\]
\end{thm}

Here $|\cdot|$ denotes the Euclidean norm. The theorem does not follow
from the usual Pr\'ekopa--Leindler inequality after subtracting
$c|x|^2/2$, since the shifted functions need not satisfy the required
pointwise hypothesis.

\medskip

In Section \ref{sec:log-sob} we turn our attention to log-Sobolev inequalities. 
The classical log-Sobolev inequality, due to Gross \cite{GrossLogarithmic1975}, may be written as
\begin{equation}
\ent(\phi)+n\log\left(\sqrt{2\pi}e\right)\int\phi
\le
\frac{1}{2}\int\frac{|\nabla\phi|^{2}}{\phi},
\label{eq:log-sobolev-orig}
\end{equation}
where
\[
\ent(\phi)=\int\phi\log\phi-\int\phi\cdot\log\int\phi
\]
is the Euclidean entropy of a function $\phi:\RR^n\to[0,\infty)$. There are
several equivalent formulations of this inequality, including the standard
one where the entropy is taken with respect to the Gaussian measure; see e.g. 
\cite[Theorem 1.3.8]{Artstein-AvidanetalAsymptotic2021}.

Bobkov and Ledoux \cite{BobkovLedouxBrunnMinkowski2000} showed that the
Pr\'ekopa--Leindler inequality also implies Gross's log-Sobolev
inequality. For log-concave functions, one obtains a more general
inequality using the same function $\alpha(t)$ from \eqref{eq:alpha-integral}. 
Since  $\log\alpha(t)$ is concave by Pr\'ekopa--Leindler, we have
\[ (\log \alpha)^\prime(0) \ge \log \alpha(1) - \log \alpha (0).\]
This gives the following theorem.

\begin{thm}\label{thm:log-sob-gen}
Let $G:\RR^n\to(-\infty,+\infty]$ be a proper lower semicontinuous
convex function such that
\[
\int_{\RR^n}e^{-G}\,dx=1.
\]
Then, for every positive integrable log-concave function
$\phi:\RR^n\to(0,\infty)$,
\[
\int_{\RR^n}
G^*\left(-\frac{\nabla\phi}{\phi}\right)\phi\,dx
\geq
n\int_{\RR^n}\phi\,dx+\ent(\phi),
\]
where
\[
G^*(y)
=
\sup_{x\in\RR^n}
\bigl(\langle x,y\rangle-G(x)\bigr)
\]
and
\[
\ent(\phi)
=
\int_{\RR^n}\phi\log\phi\,dx
-
\left(\int_{\RR^n}\phi\,dx\right)
\log\left(\int_{\RR^n}\phi\,dx\right).
\]
The gradient is understood almost everywhere, and the integral on the
left is understood in the extended sense.
\end{thm}

Taking
\[
G(x)=\frac{|x|^{2}}{2}+\frac{n}{2}\log(2\pi),
\]
we have
\[
G^*(y)=\frac{|y|^{2}}{2}-\frac{n}{2}\log(2\pi),
\]
and Theorem \ref{thm:log-sob-gen} implies the classical log-Sobolev
inequality \eqref{eq:log-sobolev-orig}, restricted here to the
log-concave case. Taking
\[
G(x)=
\begin{cases}
\log\vol(B_2^n), & |x|\le1,\\
+\infty, & |x|>1,
\end{cases}
\]
we have $G^*(y)=|y|-\log\vol(B_2^n)$, and hence obtain
\[
\ent(\phi)+\left(n+\log\vol(B_2^n)\right)\int\phi
\le \int|\nabla\phi|.
\]
Bobkov and Ledoux organized the computation in
\cite{BobkovLedouxBrunnMinkowski2000} somewhat differently, obtaining
less general inequalities but without the restriction to log-concave
functions.

In the remaining sections of this paper, we present several extensions
and variants of the Brascamp--Lieb inequality. Our interest comes from
the connection of these inequalities to convex geometry, and in particular
to the B-conjecture. In \cite{CEFM2004}, Cordero-Erausquin, Fradelizi
and Maurey proved that if $\gamma$ is the standard Gaussian measure on
$\RR^n$ and $K\subseteq\RR^n$ is an origin-symmetric convex body, then
\begin{equation}
t\mapsto\gamma(e^tK)
\label{eq:B-gauss}
\end{equation}
is log-concave. This answered a question which Lata\l{}a attributes to
Banaszczyk; see \cite{Latala2002}. Interestingly, the same result was 
proved by Brascamp and Lieb in \cite{BrascampLieb1975}, 
decades before  Lata\l{}a even suggested the 
problem, but this result remained largely unnoticed until lately -- 
see the short discussion after Theorem \ref{thm:BL-rotation}. 
In any case, the authors of \cite{CEFM2004}
conjectured that \eqref{eq:B-gauss} remains true with the Gaussian
measure replaced by any even log-concave measure on $\RR^n$. This
conjecture became known as the B-conjecture.

Very recently, after this manuscript was mostly complete, the B-conjecture 
was presumably resolved by OpenAI using an internal model. More precisely, 
their preprint \cite{openAIlogBM2026} claims a proof of the celebrated 
log-Brunn--Minkowski conjecture, originally raised and solved in dimension $n=2$
by B\"or\"oczky, Lutwak, Yang and Zhang \cite{BLYZ2012}. A result of Saroglou  
\cite{Saroglou2016} shows that the log-Brunn--Minkowski conjecture implies the
B-conjecture in full generality.  We would also like to mention the paper \cite{ENT2018}, 
where the B-conjecture was previously resolved when the measure is a certain
Gaussian mixture.  

While this recent result is a (yet to be peer reviewed) major achievement, it is
actually not very relevant to the functional inequalities discussed in this paper
and does not make them redundant in any way. Much
more relevant to us is the work \cite{CERotem2023}, where the B-conjecture was
proved for rotation-invariant log-concave measures. The proof relied
crucially on an improved Brascamp--Lieb inequality for even functions.
Let $V(x)=v(|x|)$ be a $C^2$ convex function on $\RR^n$ which depends
only on $|x|$. Let $\mu$ be the probability measure with density
$Z^{-1}e^{-V}$, and let $\nu$ be any even probability measure which is
more log-concave than $\mu$, meaning that $\frac{\dd\nu}{\dd\mu}$ is an
even log-concave function. Then, for every even $C^1$ function
$f:\RR^n\to\RR$,
\begin{equation}
\var_{\nu}f
\le
\int
\left\langle
\left(\nabla^{2}V+
\frac{v'(|x|)}{|x|}\id\right)^{-1}
\nabla f,\nabla f
\right\rangle
\dd\nu .
\label{eq:CR-BL}
\end{equation}

Since $v'(|x|)/|x|\ge0$, this is an improvement over the classical
Brascamp--Lieb inequality in the case $\nu=\mu$ and $f$ is even. For the
Gaussian measure, $v(t)=t^2/2$, and \eqref{eq:CR-BL} reduces to a
Poincar\'e inequality with constant $1/2$, which was precisely the
original result of \cite{CEFM2004}.

Using standard techniques, \eqref{eq:CR-BL} was derived in
\cite{CERotem2023} from the following Poincar\'e-type inequality: for
every $C^2$ even function $u:\RR^n\to\RR$,
\begin{equation}
\int
\frac{v'(|x|)}{|x|}|\nabla u|^2\dd\nu
\le
\int \nn{\nabla^2u}_{HS}^{2}\dd\nu ,
\label{eq:CR-poincare}
\end{equation}
where $\nn{\nabla^2u}_{HS}^{2}=\sum_{i,j}u_{x_ix_j}^{2}$ is the
Hilbert--Schmidt norm. While \cite{CERotem2023} treats the
rotation-invariant case, the same argument gives the following abstract
criterion.

\begin{proposition}\label{prop:poincare-to-B}
Let $\mu$ be an even log-concave probability measure with $C^2$ density
$Z^{-1}e^{-V}$. Assume that there exists a matrix-valued function
$T:\RR^n\to M_{n\times n}(\RR)$ with the following properties:
\begin{enumerate}
\item $T(x)$ is symmetric and $\nabla^2V(x)+T(x)$ is positive definite
for all $x\in\RR^n$.
\item $T(x)x=\nabla V(x)$ for all $x\in\RR^n$.
\item For every probability measure $\nu$ which is even and more log-concave than
$\mu$, and for every $C^2$ even function $u:\RR^n\to\RR$, one has
\begin{equation}
\int\langle T\nabla u,\nabla u\rangle\dd\nu
\le
\int\nn{\nabla^2u}_{HS}^{2}\dd\nu .
\label{eq:genT-Poincare}
\end{equation}
\end{enumerate}
Then, for all such $\nu$ and all $C^1$ even functions $f:\RR^n\to\RR$,
\begin{equation}
\var_{\nu}f
\le
\int
\left\langle
\left(\nabla^2V+T\right)^{-1}\nabla f,\nabla f
\right\rangle
\dd\nu ,
\label{eq:genT-BL}
\end{equation}
and the measure $\mu$ satisfies the B-conjecture.
\end{proposition}

We prove Proposition \ref{prop:poincare-to-B} in Section \ref{sec:unconditional}. 
This proposition illustrates how second-order improvements of the
Brascamp--Lieb inequality can have geometric consequences. While the fact that $\mu$
satisfies the B-conjecture may now be redundant by \cite{openAIlogBM2026}, we emphasize 
that the validity  of (\ref{eq:genT-BL}) is a stronger analytic fact that is 
\emph{not} established in that work, and the study of  improved Brascamp--Lieb-type 
inequalities is still very much an interesting area of active research. 

In this spirit, all the results below are motivated by this connection between functional 
inequalities and convex geometry, but the particular extensions proved here are mostly 
of independent analytic interest and are therefore unaffected
by the recent AI breakthrough. 

\begin{remark}
Before presenting our new inequalities, let us make a comment on the inequality 
\eqref{eq:CR-poincare}. 
By inspecting the proof in \cite{CERotem2023}, one sees that something stronger is actually true:
For every smooth \emph{odd} function $h:\RR^n \to \RR$ one has
\begin{equation} \label{eq:CR-poincare-improved}
\int \frac{v'(|x|)}{|x|}h^2\dd\nu + 
\int \frac{\left(\langle\nabla h, x \rangle - h\right)^2}{|x|^2} \dd \nu 
\le
\int |\nabla h|^{2}\dd\nu. 
\end{equation}
The inequality \eqref{eq:CR-poincare} then follows by neglecting the middle term, choosing $h=u_{x_i}$ and 
summing over $1 \le i \le n$. 

Assuming that $\frac{\dd \nu}{\dd \mu} = e^{-W}$ where $W$ is finite and smooth, 
the inequality \eqref{eq:CR-poincare-improved} admits the following alternative formulation: 
By performing the change of functions $h(x) = |x|f(x)$, computing both sides and integrating by
parts one obtains the inequality
\begin{equation*}
(n-1) \int f^2 \dd\nu - \int \langle \nabla W, x \rangle f^2 \dd \nu 
\le 
\int |x|^2 |\nabla_S f|^2 \dd \nu.
\end{equation*}
Here $f:\RR^n \to \RR$ is a smooth \emph{odd} function, and
\[
\nabla_S f(x) = 
\nabla f(x) - \left\langle \nabla f(x), \frac{x}{|x|} \right\rangle  \frac{x}{|x|}
\]
denotes the spherical gradient. Note that in this equivalent formulation the inequality 
no longer depends explicitly on $v$, but does depend explicitly on $W$. 
\end{remark}

In Section \ref{sec:concave}, we prove a Brascamp--Lieb inequality with a
concave weight.

\begin{thm}\label{thm:BL-concave}
Let $K\subset\RR^n$ be a convex domain, and let $V: \RR^n \to \RR$ be a $C^2$ function such
that 
\[
0<Z:=\int_{\RR^n}e^{-V}\,dx<\infty,
\qquad
\nabla^2V(x)\succeq 0\quad\text{for every }x\in K.
\]
Let
\[
d\mu=Z^{-1}e^{-V}\,dx,
\]
and let $g:K\to[0,\infty)$ be a finite concave function. Then every
$f\in C^1(K)$ satisfying
\[
\int_K |f|\,d\mu<\infty,
\qquad
\int_K f\,d\mu=0,
\]
satisfies
\begin{equation}\label{eq:BL-weighted}
\int_K gf^2\,d\mu
\leq
\int_K g
\left\langle
(\nabla^2V)^{-1}\nabla f,\nabla f
\right\rangle d\mu,
\end{equation}
where the inequality is understood in the extended sense if the
right-hand side is infinite.
\end{thm}

The fact that one may restrict the Brascamp--Lieb inequality to a convex set restriction 
to a convex set is standard. The new point in Theorem \ref{thm:BL-concave} is the 
presence of the
concave weight $g$. Since every non-negative concave function is
log-concave, one could apply Brascamp--Lieb to $V-\log g$, but this gives
the result under the centering condition $\int_K gf\dd\mu=0$. The theorem
shows that the same conclusion holds under the different condition
$\int_K f\dd\mu=0$.

In Section \ref{sec:unconditional}, we turn to the unconditional case.
A function $u:\RR^n\to\RR$ is unconditional if
\[
u(x_1,\ldots,x_n)
=
u(\epsilon_1x_1,\ldots,\epsilon_nx_n)
\]
for all $x\in\RR^n$ and all
$(\epsilon_1,\ldots,\epsilon_n)\in\{\pm1\}^n$. Unconditional measures and
sets are defined similarly. The inequality \eqref{eq:genT-Poincare}
holds for unconditional functions in a particularly strong form.

\begin{thm}\label{thm:poincare-uncond}
Let $\mu$ be an even log-concave probability measure on $\RR^n$ with $C^2$
density $Z^{-1}e^{-V}$. Fix any matrix-valued function
$T:\RR^n\to M_{n\times n}(\RR)$ such that $T(x)$ is symmetric,
$T_{ij}(x)x_ix_j\ge0$ for all $1\le i,j\le n$, and $T(x)x=\nabla V(x)$.
Let $\nu$ be any measure such that $\frac{\dd\nu}{\dd\mu}$ is
unconditional and log-concave, and let $u:\RR^n\to\RR$ be a $C^2$
unconditional function. Then
\[
\int\langle T\nabla u,\nabla u\rangle\dd\nu
\le
\int\nn{\nabla^2u}_{HS}^{2}\dd\nu .
\]
\end{thm}

If $\mu$ itself is also unconditional and $\nabla^2 V + T \succ 0$, the argument of Proposition
\ref{prop:poincare-to-B} applies on the subspace of unconditional
functions and gives the improved Brascamp--Lieb inequality
\[
\var_{\nu}f
\le
\int
\left\langle
\left(\nabla^2V+T\right)^{-1}\nabla f,\nabla f
\right\rangle\dd\nu
\]
for every admissible $T$ and every unconditional function $f$. In
particular, the map $t\mapsto\mu(e^tK)$ is log-concave when both $\mu$
and $K$ are unconditional, but this geometric consequence is not substantial: Even before
the recent AI breakthroughs, the B-conjecture in the unconditional setting had 
much simpler proofs, as was observed in \cite{CEFM2004}.

In Section \ref{sec:rotation}, we prove a general improvement of the
Brascamp--Lieb inequality for even functions.

\begin{thm}\label{thm:BL-rotation}
Let $\mu$ be an even probability measure on $\RR^n$ with density
$Z^{-1}e^{-V}$, where $V:\RR^n\to\RR$ is $C^2$ and convex. Define
\[
H_V(x)
=
\int
\left(\nabla^2V(x)+\nabla^2V(y)\right)^{-1}
\dd\mu(y).
\]
Then every $C^1$ even function $f:\RR^n\to\RR$ satisfies
\[
\var_{\mu}f
\le
\int\left\langle H_V\nabla f,\nabla f\right\rangle\dd\mu .
\]
\end{thm}

If $V$ is not strongly convex, the quadratic form defined by $H_V$ is still well-defined by
\[
\langle H_V(x)z,z\rangle = \lim_{\epsilon\to0^+}
\int \left\langle \left(\nabla^2V(x)+\nabla^2V(y)+\epsilon\id\right)^{-1}z,z
\right\rangle \dd\mu(y)
\in[0,\infty].
\]

Since $H_V\preceq(\nabla^2V)^{-1}$, Theorem \ref{thm:BL-rotation} is a strengthening of the
Brascamp--Lieb inequality. If $\nabla^2V\succeq\id$, then
$H_V\preceq\frac12\id$, and we recover the improved Poincar\'e inequality
from \cite{CEFM2004},
\[
\var_{\mu}f\le \frac12\int|\nabla f|^2\dd\mu .
\]
Surprisingly, as we already mentioned, this inequality and its connection to the B-conjecture
were known already to Brascamp and Lieb. It appeared in \cite{BrascampLieb1975}, but was not noticed by the 
community until very recently (see, e.g. the discussion in \cite{CorderoEskenazis2025}). 
The proof of \cite{BrascampLieb1975} is very different than the known one from \cite{CEFM2004}, 
and relies on a ``rotation trick''. A variant of the same trick proves Theorem \ref{thm:BL-rotation}.

Finally, in Section \ref{sec:norms}, we extend \eqref{eq:CR-poincare} to
non-Euclidean norms.

\begin{thm}\label{thm:Poincare-norms}
Let $\nn{\cdot}$ be a norm on $\RR^n$ whose unit ball
\[
K=\{x\in\RR^n:\ \nn{x}\le1\}
\]
is smooth and has strictly positive Gaussian curvature at every boundary
point. Let $v:[0,\infty)\to\RR$ be $C^1$ and increasing, and let $\mu$ be
the measure on $\RR^n$ with density $e^{-v(\nn{x})}$. Let $\nu$ be any
even measure which is more log-concave than $\mu$. Then for every odd
$C^1$ function $f:\RR^n\to\RR$,
\begin{equation}
\int\frac{v'(\nn{x})}{\nn{x}}f^2\dd\nu
\le
\int
\left\langle
\left(\nabla^2\alpha\right)^{-1}\nabla f,\nabla f
\right\rangle
\dd\nu ,
\label{eq:norm-BL}
\end{equation}
where $\alpha(x)=\frac12\nn{x}^{2}$.
\end{thm}

When $v(t)=t^2/2$, the result follows immediately from the classical
Brascamp--Lieb inequality. When $K$ is the Euclidean ball, applying
\eqref{eq:norm-BL} to $f=u_{x_i}$ and summing over $i$ recovers
\eqref{eq:CR-poincare}. For a general norm, the same argument gives, for
even $C^2$ functions $u:\RR^n\to\RR$,
\begin{equation}
\int\frac{v'(\nn{x})}{\nn{x}}|\nabla u|^2\dd\nu
\le
\int\sum_{i=1}^n
\left\langle
\left(\nabla^2\alpha\right)^{-1}\nabla u_{x_i},
\nabla u_{x_i}
\right\rangle
\dd\nu .
\label{eq:norm-order2}
\end{equation}
The right-hand side of \eqref{eq:norm-order2} does not seem to have a
simple invariant form, and we therefore do not know how to deduce an
analogue of \eqref{eq:CR-BL} from it.

\textbf{Acknowledgement.} The authors are supported by NSF-BSF
DMS-2247834 and the equivalent BSF grant 2022707. The first named author is also supported by Israel Science Foundation grant 
2574/24. The second named author is grateful to the Lady Davis Fellowship at the Technion.

\textbf{AI use statement.} The mathematical content of this paper and the initial draft were generated
solely by the authors over the course of several years, with the exception of Remarks 5.2 and 5.3. Recently, ChatGPT-5.6 pro (Sol Ultra) was used during the editorial process to identify mathematical and typographical errors, improve the writing, make some arguments more detailed, and to generate the Remarks 5.2 and 5.3. For these remarks, the authors prompted ChatGPT to provide potential applications of the results in Section 5, and amid many ideas that it suggested, the content of these two remarks was selected by the authors for inclusion in the paper.

\section{Preliminaries}\label{sec:preliminaries}

We now collect the analytic and convex-geometric facts used throughout the
paper. All integrations by parts below are first justified for smooth
compactly supported functions and smooth uniformly convex potentials. The
stated forms then follow by standard approximation arguments (see e.g. \cite[Theorem 4]{CERotem2023}). When a Hessian is only non-negative, all expressions such as 
$\langle \left(\nabla^{2}V\right)^{-1}\nabla f,\nabla f \rangle$  
are understood as limits after replacing $V$ by $V+\e |x|^2/2$ and then letting $\e\to0$.

If $\mu$ has density $Z^{-1}e^{-V}$, write
\[
L_\mu u=\Delta u-\langle\nabla V,\nabla u\rangle,
\]
so we have the integration by parts formula
\begin{equation}
\int g L_\mu u\,\dd\mu
=
-\int\langle\nabla g,\nabla u\rangle\,\dd\mu .
\label{eq:ibp}
\end{equation}
We will also use the Bochner identity
\begin{equation}
\int (L_\mu u)^2\,\dd\mu
=
\int
\left(
\nn{\nabla^2u}_{HS}^{2}
+\langle\nabla^2V\nabla u,\nabla u\rangle
\right)\dd\mu .
\label{eq:bochner}
\end{equation}

For completeness, we recall a standard proof of the Brascamp--Lieb inequality:
\begin{theorem}[Brascamp--Lieb inequality]\label{lem:BL}
Let $\mu$ be a probability measure with density $Z^{-1}e^{-V}$, where
$V\in C^2$ and $\nabla^2V\succ0$. Then for every smooth $f: \RR^n \to \RR$ we have
\[
\var_{\mu}f
\le
\int
\left\langle
(\nabla^2V)^{-1}\nabla f,\nabla f
\right\rangle\dd\mu .
\]
\end{theorem}

\begin{proof}
As already stated, we may assume by standard approximation arguments that $f$ is compactly supported. 
We may also assume without loss of generality that $\int f \dd \mu = 0$. By the standard theory
of elliptic PDEs (see e.g. \cite{Evans}) we may then find $u:\RR^n \to \RR$ smooth such that $L_\mu u=f$.

Using \eqref{eq:ibp} and the Cauchy--Schwarz inequality we obtain
\begin{align*}
\int f^2\,\dd\mu
&=
\int f L_\mu u\,\dd \mu
=
-\int\langle\nabla f,\nabla u\rangle\,\dd\mu \\
&\le
\left(
\int
\left\langle
(\nabla^2V)^{-1}\nabla f,\nabla f
\right\rangle\dd\mu
\right)^{1/2}
\left(
\int\langle\nabla^2V\nabla u,\nabla u\rangle\dd\mu
\right)^{1/2} \\
&\le
\left(
\int
\left\langle
(\nabla^2V)^{-1}\nabla f,\nabla f
\right\rangle\dd\mu
\right)^{1/2}
\left(\int f^2\,\dd\mu\right)^{1/2},
\end{align*}
where in the last inequality we also used Bochner's identity and the fact that $f=L_\mu u$. The claim now
follows immediately. 
\end{proof}

\subsection*{Legendre transform}

Given a convex function $f:\RR^n\to(-\infty,\infty]$, its Legendre
transform is
\[
f^*(x)=\sup_{y\in\RR^n}\left(\langle x,y\rangle-f(y)\right).
\]
If $f$ is smooth and strictly convex, then
\begin{equation}
f(x)+f^*(\nabla f(x))=\langle\nabla f(x),x\rangle,
\label{eq:legendre-dual}
\end{equation}
and $\nabla f^*(\nabla f(x))=x$.

Given convex functions $f_0, f_1 : \RR^n \to (-\infty,\infty])$ and $0\le t \le 1$, consider the inf-convolution $f_t$ defined by 
\[
f_t(z) = \left((1-t) f_0\star t f_1\right)(z)
=
\inf_{(1-t) x+t y=z}
\left((1-t) f_0(x)+t f_1(y)\right).
\]
Then
\begin{equation}
f_t^* = (1-t) f_0^*+ t f_1^*.
\label{eq:leg-inf-conv}
\end{equation}
In other words, the Legendre transform linearizes the inf-convolution. 

Of key importance for us will be the first and second variation formulas for the Legendre transform.
These are well-known, but we provide a proof sketch for convenience:

\begin{lemma}
\label{lem:legendre-variation}
Let \(I\subset\mathbb R\) be an open interval, and let
\[
(x,t)\longmapsto f_t(x)
\]
be a \(C^2\) function on \(\mathbb R^n\times I\). Assume that there
exists \(m>0\) such that
\[
\nabla^2 f_t(x)\succeq m\id
\]
for every \(x\in\mathbb R^n\) and \(t\in I\). Then \(f_t^*\) is finite
and \(C^2\) on \(\mathbb R^n\). Moreover, for every \(y\in\mathbb R^n\),
\[
\frac{\dd}{\dd t}f_t^*(y)
=
-\dot f_t\bigl(\nabla f_t^*(y)\bigr),
\]
and
\begin{align*}
\frac{\dd^2}{\dd t^2}f_t^*(y)
&=
-\ddot f_t\bigl(\nabla f_t^*(y)\bigr) \\
&\quad+
\left\langle
\left(
\nabla^2 f_t\bigl(\nabla f_t^*(y)\bigr)
\right)^{-1}
\nabla\dot f_t\bigl(\nabla f_t^*(y)\bigr),
\nabla\dot f_t\bigl(\nabla f_t^*(y)\bigr)
\right\rangle .
\end{align*}
\end{lemma}

\begin{proof}
The lower Hessian bound implies that \(f_t\) is superlinear and
strictly convex. Hence, for every \(y\in\mathbb R^n\), the supremum
defining \(f_t^*(y)\) is attained at a unique point
\[
x_t=\nabla f_t^*(y),
\]
characterized by
\[
\nabla f_t(x_t)=y.
\]
Thus
\[
f_t^*(y)=\langle y,x_t\rangle-f_t(x_t).
\]
Differentiating with respect to \(t\), and using
\(\nabla f_t(x_t)=y\), gives
\[
\frac{\dd}{\dd t}f_t^*(y)
=
-\dot f_t(x_t).
\]

Differentiating the identity \(\nabla f_t(x_t)=y\), we obtain
\[
\nabla\dot f_t(x_t)+\nabla^2f_t(x_t)\dot x_t=0,
\]
and therefore
\[
\dot x_t
=
-\left(\nabla^2f_t(x_t)\right)^{-1}
\nabla\dot f_t(x_t).
\]
Differentiating the first-variation formula once more now gives
\begin{align*}
\frac{\dd^2}{\dd t^2}f_t^*(y)
&=
-\ddot f_t(x_t)
-\left\langle\nabla\dot f_t(x_t),\dot x_t\right\rangle \\
&=
-\ddot f_t(x_t)
+
\left\langle
\left(\nabla^2f_t(x_t)\right)^{-1}
\nabla\dot f_t(x_t),
\nabla\dot f_t(x_t)
\right\rangle .
\end{align*}
Substituting \(x_t=\nabla f_t^*(y)\) proves the result.
\end{proof}

\section{Second variation of Pr\'ekopa--Leindler}\label{sec:PL}
In this section we prove Theorem \ref{thm:even-PL}. We begin with an improved Brascamp--Lieb inequality for 
even functions. 
\begin{proposition}\label{prop:even-BL}
Let $V:\RR^n\to\RR$ be an even $C^2$ function which satisfies  $\nabla^2V\succeq c\id$ for some $c>0$. 
Let $\mu$ be the probability measure with density $\frac{1}{Z}e^{-V}$. Then for every smooth even
function $f: \RR^n \to \RR$ we have 
\[
\var_{\mu}f
\le
\int
\left\langle
\left(\nabla^2V+c\id\right)^{-1}\nabla f,\nabla f
\right\rangle\dd\mu .
\]
\end{proposition}

Note that one can obtain Proposition \ref{prop:even-BL} as a corollary of Theorem \ref{thm:BL-rotation}.
Still, we will give a more standard $L^2$-type proof.
\begin{proof}
We may assume $f$ is compactly supported and $\int f \dd \mu = 0$. Choose $u:\RR^n \to \RR$ 
which solves $L_\mu u=f$. Since  $V$ and $f$ are even, the solution $u$ is also even. Hence each
derivative $u_i$ is odd and has mean zero. Using the fact that $\nabla^2V\succeq c\id$ and the usual 
Brascamp--Lieb inequality for $\mu$ we obtain
\[
c\int|\nabla u|^2\dd\mu
=
c\sum_{i=1}^n\int u_i^2\dd\mu
\le
\sum_{i=1}^n\int|\nabla u_i|^2\dd\mu
=
\int\nn{\nabla^2u}_{HS}^{2}\dd\mu .
\]
Using this and the Bochner identity \eqref{eq:bochner} we deduce that
\[
\int f^2\dd\mu \ge 
\int \left\langle
(\nabla^2 V +c\id)\nabla u,\nabla u
\right\rangle\dd\mu .
\]
The rest of the proof is identical to the proof given in Theorem \ref{lem:BL}. 
Using integration by parts and the Cauchy--Schwarz inequality we have
\begin{align*}
\int f^2\dd\mu &= 
\int f L_\mu u\,\dd \mu = 
-\int\langle\nabla f,\nabla u\rangle\dd\mu \\
&\le
\left(
\int
\left\langle
(\nabla^2V+c\id)^{-1}\nabla f,\nabla f
\right\rangle\dd\mu
\right)^{1/2} \\
&\quad\times
\left(
\int
\left\langle
(\nabla^2V+c\id)\nabla u,\nabla u
\right\rangle\dd\mu
\right)^{1/2} \\
&\le
\left(
\int
\left\langle
(\nabla^2V+c\id)^{-1}\nabla f,\nabla f
\right\rangle\dd\mu
\right)^{1/2}
\left(\int f^2\dd\mu\right)^{1/2}.
\end{align*}
This proves the result.
\end{proof}

Next, we explain the well-known connection between the Pr\'ekopa--Leindler inequality and 
the Brascamp--Lieb inequality. 

\begin{proposition}
\label{prop:second-variation}
Let \(f_0,f_1\in C^2(\mathbb R^n)\), and assume that there exist
constants \(0<m\leq M<\infty\) such that
\[
m\id
\preceq
\nabla^2f_0(x),\nabla^2f_1(x)
\preceq
M\id
\qquad
\text{for every }x\in\mathbb R^n.
\]
For \(0\leq t\leq1\), define
\[
f_t^*=(1-t)f_0^*+tf_1^*
\]
and
\[
\psi_t(x)
=
\bigl(f_1^*-f_0^*\bigr)\bigl(\nabla f_t(x)\bigr).
\]

Let \(\nu\) be a nonzero Borel measure on \(\mathbb R^n\). Assume that,
for some \(a\) satisfying \(0<a<m/2\),
\[
\int_{\mathbb R^n}e^{-a|x|^2}\,\dd\nu(x)<\infty.
\]
Set
\[
\alpha(t)=\int_{\mathbb R^n}e^{-f_t}\,\dd\nu
\]
and let \(\mu_t\) be the probability measure defined by
\[
\frac{\dd\mu_t}{\dd\nu}
=
\frac{e^{-f_t}}{\alpha(t)}.
\]
Then \(\alpha\) is positive and twice differentiable on \((0,1)\), and
\[
(\log\alpha)''(t)\leq0
\]
is equivalent to
\[
\var_{\mu_t}(\psi_t)
\leq
\int_{\mathbb R^n}
\left\langle
(\nabla^2f_t)^{-1}\nabla\psi_t,
\nabla\psi_t
\right\rangle
\dd\mu_t.
\]
Consequently, \(\alpha\) is log-concave on \([0,1]\) if and only if the
preceding inequality holds for every \(0<t<1\).
\end{proposition}

\begin{proof}
The Hessian bounds on \(f_0\) and \(f_1\) imply
\[
\frac1M\id
\preceq
\nabla^2f_0^*,\nabla^2f_1^*
\preceq
\frac1m\id.
\]
Therefore
\[
\frac1M\id
\preceq
\nabla^2\bigl((1-t)f_0^*+tf_1^*\bigr)
\preceq
\frac1m\id,
\]
and hence
\[
m\id\preceq\nabla^2f_t\preceq M\id
\]
for every \(0\leq t\leq1\).

Applying Lemma \ref{lem:legendre-variation} to the family
\[
(1-t)f_0^*+tf_1^*,
\]
whose Legendre transform is \(f_t\), gives
\[
\dot f_t=-\psi_t.
\]
The second-variation formula gives
\[
\ddot f_t
=
\left\langle
(\nabla^2f_t)^{-1}\nabla\psi_t,
\nabla\psi_t
\right\rangle .
\]

The uniform Hessian bounds also imply, uniformly in \(t\in[0,1]\),
\[
f_t(x)\geq \frac m2|x|^2-C(1+|x|),
\]
as well as
\[
|\psi_t(x)|\leq C(1+|x|^2),
\qquad
|\nabla\psi_t(x)|\leq C(1+|x|).
\]
Together with the assumption on \(\nu\), these estimates justify
differentiation under the integral sign. We therefore have
\[
\alpha'(t)
=
-\int_{\mathbb R^n}\dot f_t e^{-f_t}\,\dd\nu
=
\int_{\mathbb R^n}\psi_t e^{-f_t}\,\dd\nu
\]
and
\begin{align*}
\alpha''(t)
&=
\int_{\mathbb R^n}
\left(\dot f_t^2-\ddot f_t\right)e^{-f_t}\,\dd\nu \\
&=
\int_{\mathbb R^n}
\left(
\psi_t^2
-
\left\langle
(\nabla^2f_t)^{-1}\nabla\psi_t,
\nabla\psi_t
\right\rangle
\right)e^{-f_t}\,\dd\nu.
\end{align*}
Dividing by \(\alpha(t)\), we obtain
\[
(\log\alpha)''(t)
=
\var_{\mu_t}(\psi_t)
-
\int_{\mathbb R^n}
\left\langle
(\nabla^2f_t)^{-1}\nabla\psi_t,
\nabla\psi_t
\right\rangle
\dd\mu_t.
\]
This proves the equivalence.
\end{proof}

\begin{proof}[Proof of Theorem \ref{thm:even-PL}]
Write $f_0 =f$ and $f_1 = g$. It is enough to prove the result when
\[
h = f_t = ((1-t)f_0^\ast + t f_1^\ast)^\ast 
\]
because the assumption of the theorem is exactly that $h\le f_t$, and hence $e^{-h}\ge e^{-f_t}$.

By approximation, we may assume $f_0$ and $f_1$ satisfy the regularity assumptions of Proposition \ref{prop:second-variation} (Since $m\geq 2c$, the integrability hypothesis in Proposition \ref{prop:second-variation} is satisfied by choosing $c/2<a<m/2$.). The assumptions of the theorem then imply that
$\nabla^2 f_0, \nabla^2 f_1 \succeq 2c \cdot \id$. Using the relation
\[
f_t^*=(1-t)f_0^*+tf_1^*,
\]
and the formula $\nabla^2f_t^*(x) = \left(\nabla^2f_t(\nabla f_t^*(x))\right)^{-1}$, it follows
that $\nabla^2 f_t \succeq 2c \cdot \id$ for all $t \in [0,1]$. 

Next we define 
\[
\alpha(t) = \int e^{-f_t(x)}e^{\frac c2|x|^2}\,\dd x
= \int e^{-f_t} \dd \nu,
\]
where $\nu$ has density $e^{\frac{c}{2}|x|^2}$. 
By Proposition \ref{prop:second-variation} we know that the inequality $(\log \alpha)''(t) \le 0$ 
is equivalent to 
\begin{equation} \label{eq:PL-from-BL}
\var_{\mu_t}(\psi_t) \le
\int \left\langle (\nabla^2 f_t)^{-1}\nabla\psi_t,\nabla\psi_t \right\rangle \dd\mu_t,
\end{equation}
where $\psi_t=(f_1^*-f_0^*)\circ\nabla f_t$ and $\frac{\dd \mu_t}{\dd \nu} = \frac{1}{Z_t} e^{-f_t}$. In other words, $\mu_t$ has density proportional to $e^{-V_t}$, where $V_t = f_t - \frac{c}{2}|x|^2$. 
But then $\nabla^2 V_t \succeq (2c-c)\cdot \id = c\id$, and $\nabla^2 f_t = \nabla^2 V_t + c\id$, 
so the validity of \eqref{eq:PL-from-BL} follows immediately from Proposition \ref{prop:even-BL} 
(note that $\psi_t$ is even since $f_0^\ast$, $f_1^\ast$ and $f_t$ are all even).

It follows that $(\log \alpha)''(t) \le 0$ for all $t\in [0,1]$, so $\alpha$ is log-concave. Therefore
\[
\alpha(t)\ge \alpha(0)^{1-t}\alpha(1)^t,
\]
which is the desired inequality.
\end{proof}

\section{A generalized log-Sobolev inequality}\label{sec:log-sob}

We start by restating Theorem~\ref{thm:log-sob-gen}. Throughout this
section, all unspecified integrals are over $\mathbb R^n$.

\begin{theorem}[Minkowski's first inequality for functions, or generalized
log-Sobolev inequality]
\label{log-sob-mink}
Let $F:\mathbb R^n\to\mathbb R$ be convex, and let
$G:\mathbb R^n\to\mathbb R\cup\{+\infty\}$ be a proper lower
semicontinuous convex function. Assume that
\[
0<Z_F:=\int_{\mathbb R^n}e^{-F}\,dx<\infty,
\qquad
0<Z_G:=\int_{\mathbb R^n}e^{-G}\,dx<\infty.
\]
Then
\begin{equation}
\label{eq:functional-Minkowski-first}
\int_{\mathbb R^n}G^*(\nabla F)e^{-F}\,dx
\geq
\int_{\mathbb R^n}F^*(\nabla F)e^{-F}\,dx
+
Z_F\log\frac{Z_G}{Z_F}.
\end{equation}
Equivalently,
\begin{equation}
\label{eq:functional-Minkowski-first-expanded}
\int_{\mathbb R^n}G^*(\nabla F)e^{-F}\,dx
\geq
nZ_F-\int_{\mathbb R^n}Fe^{-F}\,dx
+
Z_F\log\frac{Z_G}{Z_F}.
\end{equation}
Here $\nabla F$ is understood almost everywhere, and the integral on
the left-hand side is understood as an extended integral in
$(-\infty,+\infty]$.

In particular, if $Z_F=Z_G$, then
\[
\int_{\mathbb R^n}G^*(\nabla F)e^{-F}\,dx
\geq
\int_{\mathbb R^n}F^*(\nabla F)e^{-F}\,dx.
\]
Thus, among the admissible convex functions $G$ satisfying $Z_G=Z_F$,
the left-hand side is minimized at $G=F$.
\end{theorem}

\begin{proof}
We divide the proof into the regular case and an approximation argument.

\medskip

\noindent
\textbf{Step 1: the regular case.}
Assume that $F,G\in C^2(\mathbb R^n)$ and that there exist constants
$0<m\leq M<\infty$ such that
\[
mI\leq\nabla^2F(x),\nabla^2G(x)\leq MI
\qquad
\text{for every }x\in\mathbb R^n.
\]
For $0\leq t\leq1$, define
\[
H_t=(1-t)F^*+tG^*,
\qquad
F_t=H_t^*,
\qquad
\alpha(t)=\int_{\mathbb R^n}e^{-F_t}\,dx.
\]
The Pr\'ekopa--Leindler inequality implies that
$t\mapsto\log\alpha(t)$ is concave. Since $F_0=F$ and $F_1=G$, we
obtain
\begin{equation}
\label{eq:logsob-concavity-step}
\frac{\alpha'(0)}{\alpha(0)}
\geq
\log\alpha(1)-\log\alpha(0)
=
\log\frac{Z_G}{Z_F}.
\end{equation}

For every $x$ and $t$, the supremum defining $F_t(x)$ has a unique
maximizer $y_t(x)$ under our assumption. Thus
\[
F_t(x)
=
\langle x,y_t(x)\rangle
-
(1-t)F^*(y_t(x))
-
tG^*(y_t(x)).
\]
Differentiating this identity in $t$, and using Lemma \ref{lem:legendre-variation}, gives
\[
\frac{\partial}{\partial t}F_t(x)
=
F^*(y_t(x))-G^*(y_t(x)).
\]
At $t=0$, the maximizing point is
\[
y_0(x)=\nabla F(x).
\]
Consequently,
\[
\left.\frac{\partial}{\partial t}\right|_{t=0+}F_t(x)
=
F^*(\nabla F(x))-G^*(\nabla F(x)).
\]

The two-sided Hessian bounds imply uniform quadratic estimates for
$F_t$, for $t$ in a neighborhood of zero. More precisely, there are
constants $c,C>0$ such that
\[
F_t(x)\geq c|x|^2-C
\]
and
\[
\left|\frac{\partial}{\partial t}F_t(x)\right|
\leq C(1+|x|^2).
\]
Therefore differentiation under the integral sign is justified, and
\[
\alpha'(0)
=
\int_{\mathbb R^n}
\bigl(G^*(\nabla F)-F^*(\nabla F)\bigr)e^{-F}\,dx.
\]
Substitution into \eqref{eq:logsob-concavity-step} gives
\[
\int_{\mathbb R^n}G^*(\nabla F)e^{-F}\,dx
\geq
\int_{\mathbb R^n}F^*(\nabla F)e^{-F}\,dx
+
Z_F\log\frac{Z_G}{Z_F}.
\]

In view of our regularity assumption, Fenchel's identity and integration by
parts give
\[
F^*(\nabla F(x))
=
\langle x,\nabla F(x)\rangle-F(x)
\]
and
\[
\int_{\mathbb R^n}
\langle x,\nabla F(x)\rangle e^{-F(x)}\,dx
=
nZ_F.
\]
Hence
\[
\int_{\mathbb R^n}F^*(\nabla F)e^{-F}\,dx
=
nZ_F-\int_{\mathbb R^n}Fe^{-F}\,dx.
\]
This proves both formulations in the regular case.

\medskip

\noindent
\textbf{Step 2: approximation.}

We first record some standard integrability facts. Put
\[
f=e^{-F}.
\]
Since $F$ is finite and convex and $f$ is integrable, there exist
constants $c,C>0$ such that
\[
f(x)\leq Ce^{-c|x|}.
\]
Moreover,
\begin{equation}
\label{eq:logconcave-W11}
f\in W^{1,1}(\mathbb R^n),
\qquad
\int_{\mathbb R^n}(1+|x|)|\nabla f(x)|\,dx<\infty.
\end{equation}
For completeness, these assertions follow from the coarea formula.
Indeed, the superlevel sets of $f$ are bounded convex sets contained
in balls whose radii are bounded by
$C(1+|\log t|)$. Monotonicity of the surface area of convex bodies
under inclusion therefore gives
\[
\int_0^{\|f\|_\infty}
(1+|\log t|)^n\,dt<\infty
\]
as an upper bound for the relevant coarea integrals.

In particular,
\begin{equation}
\label{eq:score-zero-mean}
\int_{\mathbb R^n}\nabla F\,e^{-F}\,dx
=
-\int_{\mathbb R^n}\nabla f\,dx
=
0.
\end{equation}
Also, $Fe^{-F}=-f\log f$ is integrable.

If $x_0\in\operatorname{dom}G$, then
\[
G^*(y)\geq\langle x_0,y\rangle-G(x_0).
\]
Together with \eqref{eq:logconcave-W11}, this shows that the negative
part of $G^*(\nabla F)e^{-F}$ is integrable. Thus the left-hand side
of the theorem is well-defined as an extended integral.

\smallskip

\noindent
\emph{Part 2a: regular $F$ and $G$ with bounded effective domain.}
Assume that $F$ satisfies the regularity assumptions of Step 1 and
that
\[
D:=\overline{\operatorname{dom}G}
\]
is compact. Set
\[
H=G^*.
\]
Then $H$ is finite and globally Lipschitz. Let $L$ be a Lipschitz
constant for $H$. Choose a smooth nonnegative symmetric mollifier
$\rho$, supported in the Euclidean unit ball and satisfying
$\int\rho=1$, and put
\[
\rho_\delta(y)=\delta^{-n}\rho(y/\delta),
\qquad
K_\delta=H\ast\rho_\delta,
\qquad
H_\delta(y)=K_\delta(y)+\frac{\delta}{2}|y|^2.
\]
Since $\rho_\delta$ is centered, convexity and Lipschitz continuity of
$H$ give
\begin{equation}
\label{eq:H-mollification-bounds}
H(y)
\leq
K_\delta(y)
\leq
H(y)+c_\rho L\delta,
\end{equation}
where $c_\rho$ depends only on $\rho$. Moreover,
\[
\delta I
\leq
\nabla^2H_\delta
\leq
\left(\delta+\frac{C_\rho L}{\delta}\right)I.
\]
Define
\[
G_\delta=H_\delta^*.
\]
The preceding two-sided Hessian bounds imply corresponding two-sided
positive Hessian bounds for $G_\delta$. Thus $G_\delta$ satisfies the
regularity assumptions of Step 1.

We claim that
\begin{equation}
\label{eq:ZG-delta-limit}
\int_{\mathbb R^n}e^{-G_\delta}\,dx
\longrightarrow
Z_G.
\end{equation}
Let
\[
G^{\mathrm{MY}}_\delta(x)
=
\inf_{z\in\mathbb R^n}
\left\{
G(z)+\frac{|x-z|^2}{2\delta}
\right\}
\]
be the Moreau envelope of $G$. From
\eqref{eq:H-mollification-bounds} and the order-reversing property of
the Legendre transform, we obtain
\[
G-c_\rho L\delta
\leq
K_\delta^*
\leq
G.
\]
Since
\[
G_\delta
=
K_\delta^*
\mathbin{\square}
\frac{|\cdot|^2}{2\delta},
\]
it follows that
\begin{equation}
\label{eq:G-delta-sandwich}
G^{\mathrm{MY}}_\delta-c_\rho L\delta
\leq
G_\delta
\leq
G.
\end{equation}
Because $G^{\mathrm{MY}}_\delta\to G$ pointwise in the extended sense,
we have $G_\delta\to G$ pointwise in the extended sense.

Since $D$ is compact and $G$ is bounded below on its effective domain,
there is $b\in\mathbb R$ such that
\[
G^{\mathrm{MY}}_\delta(x)
\geq
b+\frac{\operatorname{dist}(x,D)^2}{2\delta}.
\]
Consequently, for $0<\delta\leq1$, the functions $e^{-G_\delta}$ are
dominated by a fixed integrable function proportional to
\[
\exp\left(-\frac{\operatorname{dist}(x,D)^2}{2}\right).
\]
Dominated convergence proves \eqref{eq:ZG-delta-limit}.

Applying Step 1 to $F$ and $G_\delta$, whose conjugate is $H_\delta$,
gives
\[
\int_{\mathbb R^n}H_\delta(\nabla F)e^{-F}\,dx
\geq
nZ_F-\int_{\mathbb R^n}Fe^{-F}\,dx
+
Z_F\log
\frac{\int e^{-G_\delta}}{Z_F}.
\]
On the other hand, \eqref{eq:H-mollification-bounds} gives
\[
H_\delta(y)
\leq
H(y)+c_\rho L\delta+\frac{\delta}{2}|y|^2.
\]
Since the Hessian bounds on $F$ imply
\[
\int_{\mathbb R^n}|\nabla F|^2e^{-F}\,dx<\infty,
\]
we may let $\delta\downarrow0$ and obtain
\begin{equation}
\label{eq:regular-F-compact-G}
\int_{\mathbb R^n}G^*(\nabla F)e^{-F}\,dx
\geq
nZ_F-\int_{\mathbb R^n}Fe^{-F}\,dx
+
Z_F\log\frac{Z_G}{Z_F}.
\end{equation}

\smallskip

\noindent
\emph{Part 2b: arbitrary $F$ and $G$ with bounded effective domain.}
Continue to assume that
$\overline{\operatorname{dom}G}$ is compact, and write again
\[
H=G^*.
\]
Thus $H$ is finite and globally Lipschitz.

Let $\gamma_s$ denote the centered Gaussian density with covariance
$sI$, and define
\[
f_s=f\ast\gamma_s,
\qquad
F_s=-\log f_s.
\]
By Pr\'ekopa's theorem, $F_s$ is convex. It is also smooth, and the
Gaussian differentiation formula gives
\[
0\leq\nabla^2F_s\leq\frac1sI.
\]
For $\varepsilon>0$, define
\[
F_{s,\varepsilon}(x)
=
F_s(x)+\frac{\varepsilon}{2}|x|^2.
\]
Then
\[
\varepsilon I
\leq
\nabla^2F_{s,\varepsilon}
\leq
\left(\frac1s+\varepsilon\right)I.
\]
Therefore Part 2a applies to $F_{s,\varepsilon}$ and $G$.

Letting $\varepsilon\downarrow0$, by monotone convergence theorem, we have
\[
\int_{\mathbb R^n}e^{-F_{s,\varepsilon}}\,dx
\longrightarrow
\int_{\mathbb R^n}f_s\,dx
=
Z_F
\]
and, by dominated convergence,
\[
\int_{\mathbb R^n}
F_{s,\varepsilon}e^{-F_{s,\varepsilon}}\,dx
\longrightarrow
\int_{\mathbb R^n}F_sf_s\,dx.
\]
Since $H$ is Lipschitz, $\nabla F_s$ has at most linear growth, and
$f_s$ has exponential tails, we also have
\[
\int_{\mathbb R^n}
H(\nabla F_{s,\varepsilon})e^{-F_{s,\varepsilon}}\,dx
\longrightarrow
\int_{\mathbb R^n}H(\nabla F_s)f_s\,dx.
\]
It follows from \eqref{eq:regular-F-compact-G} that
\begin{equation}
\label{eq:regularized-F-inequality}
\int_{\mathbb R^n}H(\nabla F_s)f_s\,dx
\geq
nZ_F-\int_{\mathbb R^n}F_sf_s\,dx
+
Z_F\log\frac{Z_G}{Z_F}.
\end{equation}

We now compare the regularized left-hand side with the original one.
Since
\[
\nabla f=-f\nabla F,
\]
differentiation of the convolution gives
\[
\nabla F_s(x)f_s(x)
=
\int_{\mathbb R^n}
\nabla F(y)f(y)\gamma_s(x-y)\,dy.
\]
Thus $\nabla F_s(x)$ is the barycenter of $\nabla F(y)$ with respect
to the probability measure
\[
\frac{f(y)\gamma_s(x-y)}{f_s(x)}\,dy.
\]
Convexity of $H$ therefore yields
\[
H(\nabla F_s(x))f_s(x)
\leq
\int_{\mathbb R^n}
H(\nabla F(y))f(y)\gamma_s(x-y)\,dy.
\]
Integrating in $x$ gives the crucial estimate
\begin{equation}
\label{eq:score-Jensen}
\int_{\mathbb R^n}H(\nabla F_s)f_s\,dx
\leq
\int_{\mathbb R^n}H(\nabla F)f\,dx.
\end{equation}

Finally, $f_s\to f$ pointwise and locally uniformly. The exponential
tail estimate for $f$ implies, uniformly for $0<s\leq1$, that
\[
f_s(x)\leq C'e^{-c'|x|}.
\]
It follows by dominated convergence that
\[
\int_{\mathbb R^n}f_s\log f_s\,dx
\longrightarrow
\int_{\mathbb R^n}f\log f\,dx,
\]
or equivalently,
\begin{equation}
\label{eq:entropy-convolution-limit}
\int_{\mathbb R^n}F_sf_s\,dx
\longrightarrow
\int_{\mathbb R^n}Fe^{-F}\,dx.
\end{equation}
Combining \eqref{eq:regularized-F-inequality},
\eqref{eq:score-Jensen}, and
\eqref{eq:entropy-convolution-limit}, and then letting $s\downarrow0$,
proves \eqref{eq:functional-Minkowski-first-expanded} whenever
$\overline{\operatorname{dom}G}$ is compact.

\smallskip

\noindent
\emph{Part 2c: arbitrary $G$.}
Choose $x_0\in\operatorname{dom}G$. For all sufficiently large $R$,
let
\[
G_R=G+\iota_{\overline{B_R}},
\]
where $\iota_{\overline{B_R}}$ is zero on $\overline{B_R}$ and
$+\infty$ outside $\overline{B_R}$. Then
\[
e^{-G_R}=e^{-G}\mathbf 1_{\overline{B_R}},
\qquad
Z_{G_R}\uparrow Z_G.
\]
Moreover,
\[
G_R^*(y)
=
\sup_{x\in\overline{B_R}}
\bigl\{\langle x,y\rangle-G(x)\bigr\}
\uparrow
G^*(y).
\]
For every $R\geq |x_0|$,
\[
G_R^*(y)
\geq
\langle x_0,y\rangle-G(x_0).
\]
Consequently,
\[
G_R^*(\nabla F)
-
\langle x_0,\nabla F\rangle
+
G(x_0)
\]
is a nonnegative increasing sequence. By monotone convergence and
\eqref{eq:logconcave-W11},
\[
\int_{\mathbb R^n}G_R^*(\nabla F)e^{-F}\,dx
\uparrow
\int_{\mathbb R^n}G^*(\nabla F)e^{-F}\,dx.
\]
Applying the already proved compact-domain case to $G_R$ and letting
$R\to\infty$ proves
\eqref{eq:functional-Minkowski-first-expanded} for arbitrary $G$.

It remains to identify the first two terms on the right-hand side.
The exponential tail estimate implies that $F$ is coercive, so it has
a minimizer $a\in\mathbb R^n$. By convexity,
\[
\langle x-a,\nabla F(x)\rangle\geq0
\]
at every differentiability point of $F$.

Let $\chi_R$ be radial cutoff functions centered at $a$, increasing
pointwise to one, such that
\[
|(x-a)\cdot\nabla\chi_R(x)|\leq C
\]
and $\nabla\chi_R$ is supported in
$\{R\leq|x-a|\leq2R\}$. Integration by parts gives
\[
\begin{aligned}
\int_{\mathbb R^n}
\chi_R\langle x-a,\nabla F\rangle e^{-F}\,dx
&=
-\int_{\mathbb R^n}
\chi_R\langle x-a,\nabla f\rangle\,dx
\\
&=
n\int_{\mathbb R^n}\chi_Rf\,dx
+
\int_{\mathbb R^n}
\langle x-a,\nabla\chi_R\rangle f\,dx.
\end{aligned}
\]
The last term tends to zero as $R\to\infty$. Since the integrand on
the left is nonnegative, monotone convergence yields
\[
\int_{\mathbb R^n}
\langle x-a,\nabla F\rangle e^{-F}\,dx
=
nZ_F.
\]
Together with \eqref{eq:score-zero-mean}, this gives
\[
\int_{\mathbb R^n}
\langle x,\nabla F\rangle e^{-F}\,dx
=
nZ_F.
\]
Finally, Fenchel's identity gives almost everywhere
\[
F^*(\nabla F(x))
=
\langle x,\nabla F(x)\rangle-F(x).
\]
Therefore
\[
\int_{\mathbb R^n}F^*(\nabla F)e^{-F}\,dx
=
nZ_F-\int_{\mathbb R^n}Fe^{-F}\,dx.
\]
Substituting this identity into
\eqref{eq:functional-Minkowski-first-expanded} proves
\eqref{eq:functional-Minkowski-first}.
\end{proof}

\begin{remark}
\label{rem:recover-Minkowski}
Let $K,L\subset\mathbb R^n$ be convex bodies. We first assume that
$|K|=|L|$. For $m>0$, set
\[
F_m=m d_K,
\qquad
G_m=m d_L,
\]
and
\[
A_m=\int_{\mathbb R^n}e^{-m d_K}\,dx,
\qquad
B_m=\int_{\mathbb R^n}e^{-m d_L}\,dx.
\]
The conjugates of these functions are
\[
F_m^*=h_K+\iota_{mB_2^n},
\qquad
G_m^*=h_L+\iota_{mB_2^n}.
\]
Since $\nabla F_m=m\nabla d_K\in mB_2^n$ almost everywhere,
Theorem~\ref{log-sob-mink} gives
\begin{equation}
\label{eq:distance-Minkowski}
\int_{\mathbb R^n}
\bigl(h_L(m\nabla d_K)-h_K(m\nabla d_K)\bigr)e^{-m d_K}\,dx
\geq
A_m\log\frac{B_m}{A_m}.
\end{equation}

For any convex body $M$, the coarea formula gives
\begin{align*}
\int_{\mathbb R^n}
h_M(m\nabla d_K)e^{-m d_K}\,dx
&=
m\int_0^\infty e^{-mr}
\int_{\partial(K+rB_2^n)}
h_M(n_{K+rB_2^n})\,d\mathcal H^{n-1}\,dr\\
&=
n\int_0^\infty e^{-s}
V_1\left(K+\frac{s}{m}B_2^n,M\right)\,ds.
\end{align*}
Here we used
\[
\int_{\partial C}h_M(n_C)\,d\mathcal H^{n-1}
=
nV_1(C,M).
\]
Therefore, by continuity of mixed volumes and dominated convergence,
\[
\lim_{m\to\infty}
\int_{\mathbb R^n}
h_M(m\nabla d_K)e^{-m d_K}\,dx
=
nV_1(K,M).
\]
Indeed, for $m\geq1$ the integrand is dominated by
\[
ne^{-s}V_1(K+sB_2^n,M),
\]
which is integrable because $V_1(K+sB_2^n,M)$ is a polynomial in $s$.

Moreover,
\[
A_m\longrightarrow |K|,
\qquad
B_m\longrightarrow |L|,
\]
by dominated convergence. Since $|K|=|L|$, the right-hand side of
\eqref{eq:distance-Minkowski} tends to zero. Passing to the limit gives
\[
V_1(K,L)\geq V_1(K,K)=|K|.
\]

For arbitrary $K$ and $L$, apply the equal-volume case to $K$ and
$\lambda L$, where
\[
\lambda=\left(\frac{|K|}{|L|}\right)^{1/n}.
\]
Using $V_1(K,\lambda L)=\lambda V_1(K,L)$, we obtain
\[
V_1(K,L)
\geq
|K|^{\frac{n-1}{n}}|L|^{\frac1n},
\]
which is Minkowski's first inequality.
\end{remark}

\medskip

Recall that, for a nonnegative integrable function $\phi$ and a
probability measure $\mu$, its entropy relative to $\mu$ is
\[
\operatorname{Ent}_{\mu}(\phi)
:=
\int\phi\log\phi\,d\mu
-\left(\int\phi\,d\mu\right)
 \log\left(\int\phi\,d\mu\right),
\]
where $0\log0=0$. Jensen's inequality gives
$\operatorname{Ent}_{\mu}(\phi)\geq0$. For Lebesgue measure we write
$\operatorname{Ent}(\phi)$. Notice that entropy relative to Lebesgue
measure need not be nonnegative.

If $Z_G=1$ and $\phi=e^{-F}$, then
\[
\operatorname{Ent}(\phi)
=
-\int_{\mathbb R^n}F e^{-F}\,dx
-\left(\int_{\mathbb R^n}e^{-F}\,dx\right)
 \log\left(\int_{\mathbb R^n}e^{-F}\,dx\right).
\]
Thus Theorem~\ref{log-sob-mink} implies the following.

\begin{cor}
\label{cor:generalized-log-sobolev-entropy}
Let $\phi=e^{-F}$ be a nonzero integrable log-concave function, and
let $G:\mathbb R^n\to\mathbb R\cup\{+\infty\}$ be a proper lower
semicontinuous convex function such that
\[
\int_{\mathbb R^n}e^{-G}\,dx=1.
\]
Under the hypotheses of Theorem~\ref{log-sob-mink},
\[
\int_{\mathbb R^n}
G^*\left(-\frac{\nabla\phi}{\phi}\right)\phi\,dx
\geq
n\int_{\mathbb R^n}\phi\,dx
+\operatorname{Ent}(\phi).
\]
The quotient $-\nabla\phi/\phi$ is taken on $\{\phi>0\}$; the
nonsmooth cases below follow by convex approximation.
\end{cor}

We now record several special choices of $G$.

\medskip

Let $\omega_n=|B_2^n|$ and choose
\[
G^*(y)=|y|-\log\omega_n.
\]
Then
\[
G(x)=\iota_{B_2^n}(x)+\log\omega_n,
\qquad
\int_{\mathbb R^n}e^{-G}\,dx=1.
\]
The preceding corollary therefore gives the following result; see also
Bobkov and Ledoux~\cite{BobkovLedouxBrunnMinkowski2000}.

\begin{cor}[$L^1$ log-Sobolev inequality]
\label{cor:L1-log-sobolev}
Let $\phi$ be a nonzero nonnegative integrable log-concave function in
$BV(\mathbb R^n)$. Then
\[
\operatorname{Ent}(\phi)
+\bigl(n+\log\omega_n\bigr)
 \int_{\mathbb R^n}\phi\,dx
\leq
|D\phi|(\mathbb R^n).
\]
In particular, if $\phi\in W^{1,1}(\mathbb R^n)$, then
\[
\operatorname{Ent}(\phi)
+\bigl(n+\log\omega_n\bigr)
 \int_{\mathbb R^n}\phi\,dx
\leq
\int_{\mathbb R^n}|\nabla\phi|\,dx.
\]
\end{cor}

\medskip

Alternatively, choose
\[
G(x)=\frac{|x|^2}{2}+\frac n2\log(2\pi).
\]
Then
\[
G^*(y)=\frac{|y|^2}{2}-\frac n2\log(2\pi),
\qquad
\int_{\mathbb R^n}e^{-G}\,dx=1.
\]
For a positive log-concave function $\phi$, the entropy form of
Theorem~\ref{log-sob-mink} gives
\begin{equation}
\label{eq:Lebesgue-log-Sobolev-phi}
\operatorname{Ent}(\phi)
+n\log\bigl(\sqrt{2\pi}\,e\bigr)
 \int_{\mathbb R^n}\phi\,dx
\leq
\frac12\int_{\mathbb R^n}
\frac{|\nabla\phi|^2}{\phi}\,dx.
\end{equation}
Writing $\phi=f^2$ gives
\begin{equation}
\label{eq:Lebesgue-log-Sobolev-f}
\operatorname{Ent}(f^2)
+n\log\bigl(\sqrt{2\pi}\,e\bigr)
 \int_{\mathbb R^n}f^2\,dx
\leq
2\int_{\mathbb R^n}|\nabla f|^2\,dx
\end{equation}
for nonnegative log-concave $f$. In fact,
\eqref{eq:Lebesgue-log-Sobolev-f} is valid for every
$f\in W^{1,2}(\mathbb R^n)$ and is equivalent to the classical
Gaussian log-Sobolev inequality (in the partial case of log-concave functions):

\begin{cor}[Classical Gaussian log-Sobolev inequality]
\label{CLLSI}
Let $\gamma$ be the standard Gaussian probability measure on
$\mathbb R^n$. Then, for every log-concave $g\in W^{1,2}(\gamma)$,
\[
\operatorname{Ent}_{\gamma}(g^2)
\leq
2\int_{\mathbb R^n}|\nabla g|^2\,d\gamma.
\]
\end{cor}
\begin{remark}
In fact the above result is classical and holds for all $W^{1,2}(\gamma)$, not just log-concave ones.
\end{remark}

Indeed, let
\[
\rho(x)=(2\pi)^{-n/2}e^{-|x|^2/2}
\]
be the Gaussian density and set $f=g\sqrt{\rho}$. A direct integration
by parts gives the identity
\begin{align*}
&2\int_{\mathbb R^n}|\nabla f|^2\,dx
-\operatorname{Ent}(f^2)
-n\log\bigl(\sqrt{2\pi}\,e\bigr)
 \int_{\mathbb R^n}f^2\,dx
\\
&\qquad=
2\int_{\mathbb R^n}|\nabla g|^2\,d\gamma
-\operatorname{Ent}_{\gamma}(g^2),
\end{align*}
which proves the equivalence.

\medskip

More generally, let $p>1$, let $q=p/(p-1)$, and set
\[
Z_{n,q}
=
\int_{\mathbb R^n}e^{-|x|^q/q}\,dx
=
|S^{n-1}|\int_0^\infty r^{n-1}e^{-r^q/q}\,dr.
\]
Choose
\[
G^*(y)=\frac{|y|^p}{p}-\log Z_{n,q}.
\]
Since
\[
G(x)=\frac{|x|^q}{q}+\log Z_{n,q},
\]
we have $\int e^{-G}=1$.

\begin{cor}[$L^p$ log-Sobolev inequality]
\label{cor:Lp-log-sobolev}
Let $p>1$, let $q=p/(p-1)$, and let $f$ be a nonzero nonnegative
log-concave function in $W^{1,p}(\mathbb R^n)$. Then
\[
\operatorname{Ent}(f^p)
+C_{n,p}\int_{\mathbb R^n}f^p\,dx
\leq
p^{p-1}\int_{\mathbb R^n}|\nabla f|^p\,dx,
\]
where
\[
C_{n,p}=n+\log Z_{n,q}.
\]
Equivalently,
\[
C_{n,p}
=
n+\log\left(
|S^{n-1}|\int_0^\infty
r^{n-1}e^{-r^q/q}\,dr
\right).
\]
\end{cor}

\begin{proof}
Apply Corollary~\ref{cor:generalized-log-sobolev-entropy} to
$\phi=f^p$. On $\{f>0\}$,
\[
-\frac{\nabla\phi}{\phi}
=
-p\frac{\nabla f}{f},
\]
and hence
\begin{align*}
\int_{\mathbb R^n}
G^*\left(-\frac{\nabla\phi}{\phi}\right)\phi\,dx
&=
p^{p-1}\int_{\mathbb R^n}|\nabla f|^p\,dx
\\
&\quad
-(\log Z_{n,q})\int_{\mathbb R^n}f^p\,dx.
\end{align*}
The conclusion follows.
\end{proof}

The same $L^p$ inequality is valid without the log-concavity
assumption, by the sharp Euclidean $L^p$ log-Sobolev inequality. The
argument above only proves it directly on the log-concave class.

\medskip

We finish with a version relative to an arbitrary log-concave
probability measure.

\begin{theorem}[Weighted generalized log-Sobolev inequality]
\label{Log-Sob-measure}
Let
\[
d\mu(x)=Z^{-1}e^{-V(x)}\,dx,
\]
where \(V:\mathbb R^n\to\mathbb R\cup\{+\infty\}\) is proper, lower
semicontinuous and convex, and \(0<Z<\infty\). Let
\(f:\mathbb R^n\to\mathbb R\) be convex and continuously
differentiable, and let
\(G:\mathbb R^n\to\mathbb R\cup\{+\infty\}\) be proper, lower
semicontinuous and convex. Assume that
\[
0<\int e^{-f}\,d\mu<\infty,
\qquad
0<\int e^{-G}\,d\mu<\infty,
\]
and that, for some \(\delta>0\),
\[
\int
\exp\left(
-f+\delta
\bigl[G^*(\nabla f)-f^*(\nabla f)\bigr]_+
\right)\,d\mu<\infty,
\]
while
\[
\int
\bigl[G^*(\nabla f)-f^*(\nabla f)\bigr]_-
e^{-f}\,d\mu<\infty.
\]
Then
\begin{align}
\int G^*(\nabla f)e^{-f}\,d\mu
\geq{}&
\operatorname{Ent}_{\mu}(e^{-f})
+\int \langle x,\nabla f(x)\rangle e^{-f}\,d\mu
\nonumber\\
&+
\left(\int e^{-f}\,d\mu\right)
\log\left(\int e^{-G}\,d\mu\right).
\label{eq:weighted-generalized-log-sobolev}
\end{align}
In particular, if
\[
\int e^{-G}\,d\mu=1
\]
and \(f\) attains its minimum at the origin, then
\[
\int G^*(\nabla f)e^{-f}\,d\mu
\geq
\operatorname{Ent}_{\mu}(e^{-f}).
\]
This applies, in particular, when \(f\) is even.
\end{theorem}

\begin{proof}
For \(0<t<1\), let
\[
f_t=\bigl((1-t)f^*+tG^*\bigr)^*.
\]
For every \(x,y\in\mathbb R^n\), the definition of the Legendre
transform gives
\[
f_t((1-t)x+ty)
\leq
(1-t)f(x)+tG(y).
\]
Combining this with the convexity of \(V\) and applying the
Pr\'ekopa--Leindler inequality, we obtain
\begin{equation}
\int e^{-f_t}\,d\mu
\geq
\left(\int e^{-f}\,d\mu\right)^{1-t}
\left(\int e^{-G}\,d\mu\right)^t.
\label{eq:weighted-PL-interpolation}
\end{equation}

On the other hand, taking \(\nabla f(x)\) as a competitor in the
supremum defining \(f_t(x)\), we find
\[
f_t(x)
\geq
f(x)-t\bigl(G^*(\nabla f(x))-f^*(\nabla f(x))\bigr).
\]
It follows that
\[
\left(\int e^{-f}\,d\mu\right)^{1-t}
\left(\int e^{-G}\,d\mu\right)^t
\leq
\int
\exp\left(
-f+t\bigl(G^*(\nabla f)-f^*(\nabla f)\bigr)
\right)\,d\mu.
\]
The two sides agree at \(t=0\). The integrability assumptions allow
us to differentiate this inequality from the right at \(t=0\), giving
\[
\left(\int e^{-f}\,d\mu\right)
\log
\frac{\int e^{-G}\,d\mu}{\int e^{-f}\,d\mu}
\leq
\int
\bigl(G^*(\nabla f)-f^*(\nabla f)\bigr)e^{-f}\,d\mu.
\]
Finally,
\[
f^*(\nabla f(x))
=
\langle x,\nabla f(x)\rangle-f(x)
\]
and
\[
\operatorname{Ent}_{\mu}(e^{-f})
=
-\int f e^{-f}\,d\mu
-
\left(\int e^{-f}\,d\mu\right)
\log\left(\int e^{-f}\,d\mu\right).
\]
Substitution gives
\eqref{eq:weighted-generalized-log-sobolev}.

If \(f\) attains its minimum at the origin, convexity gives
\[
\langle x,\nabla f(x)\rangle
\geq f(x)-f(0)\geq0.
\]
Thus, when \(\int e^{-G}\,d\mu=1\), the additional terms on the
right-hand side of
\eqref{eq:weighted-generalized-log-sobolev} are nonnegative.
Finally, every even convex function attains its minimum at the origin.
\end{proof}

Related inequalities obtained by Pr\'ekopa--Leindler methods were
studied by Bobkov and Ledoux
\cite{BobkovLedouxBrunnMinkowski2000}.

\section{A Brascamp--Lieb inequality with a concave weight}
\label{sec:concave}

We begin with a weighted Bochner identity (see \cite{KolesnikovMilman2018}, \cite{KolLiv}).

\begin{lemma}[Weighted Bochner identity]
\label{lem:weighted-bochner}
Let $\dd\mu=e^{-V}\dd x$, where $V\in C^2(\mathbb R^n)$, and let
$g,u$ be smooth. Assume that all the terms below are integrable and that
the integrations by parts have no boundary contribution. Then
\begin{align*}
\int_{\mathbb R^n} g(L_\mu u)^2\dd\mu
&=
\int_{\mathbb R^n} g
\left(
\nn{\nabla^2u}_{HS}^{2}
+\left\langle\nabla^2V\nabla u,\nabla u\right\rangle
\right)\dd\mu \\
&\quad
-\int_{\mathbb R^n}
\left\langle\nabla^2g\nabla u,\nabla u\right\rangle\dd\mu
-2\int_{\mathbb R^n}
\left\langle\nabla g,\nabla u\right\rangle L_\mu u\,\dd\mu .
\end{align*}
\end{lemma}

\begin{proof}[Proof for completeness]
Recall that
\[
L_\mu=\Delta-\left\langle\nabla V,\nabla\right\rangle
\]
and that
\[
\frac12L_\mu|\nabla u|^2
=
\left\langle\nabla u,\nabla L_\mu u\right\rangle
+\nn{\nabla^2u}_{HS}^{2}
+\left\langle\nabla^2V\nabla u,\nabla u\right\rangle .
\]
Integration by parts gives
\begin{align*}
\int g(L_\mu u)^2\dd\mu
&=
-\int L_\mu u\,
\left\langle\nabla g,\nabla u\right\rangle\dd\mu
-\int g\left\langle\nabla L_\mu u,\nabla u\right\rangle\dd\mu \\
&=
-\int L_\mu u\,
\left\langle\nabla g,\nabla u\right\rangle\dd\mu \\
&\quad
+\int g
\left(
\nn{\nabla^2u}_{HS}^{2}
+\left\langle\nabla^2V\nabla u,\nabla u\right\rangle
\right)\dd\mu
-\frac12\int gL_\mu|\nabla u|^2\dd\mu .
\end{align*}
Moreover,
\[
-\frac12\int gL_\mu|\nabla u|^2\dd\mu
=
\frac12\int
\left\langle\nabla g,\nabla|\nabla u|^2\right\rangle\dd\mu
=
\int
\left\langle\nabla^2u\nabla g,\nabla u\right\rangle\dd\mu .
\]
Finally,
\begin{align*}
0
&=
\int
\div_\mu\left(
\left\langle\nabla g,\nabla u\right\rangle\nabla u
\right)\dd\mu \\
&=
\int
\left(
\left\langle\nabla^2g\nabla u,\nabla u\right\rangle
+\left\langle\nabla^2u\nabla g,\nabla u\right\rangle
+\left\langle\nabla g,\nabla u\right\rangle L_\mu u
\right)\dd\mu .
\end{align*}
Substitution gives the asserted identity.
\end{proof}

\begin{proof}[Proof of Theorem \ref{thm:BL-concave}]
We first prove the assertion on a bounded convex domain $D$ with smooth
boundary. Assume that $V$ and $f$ are smooth on a neighborhood of
$\overline D$, that $g\in C^2(\overline D)$ is nonnegative and concave,
and that
\[
\int_D f\,\dd\mu=0.
\]
Let $u$ be the solution, unique up to an additive constant, of the
weighted Neumann problem
\[
\begin{cases}
L_\mu u=f & \text{in }D,\\
\partial_\nu u=0 & \text{on }\partial D,
\end{cases}
\]
where $\nu$ is the outward unit normal. The displayed centering
condition is precisely the compatibility condition for this problem.

We use the convention
\[
\mathrm{II}(X,Y)=\left\langle D_X\nu,Y\right\rangle
\]
for the second fundamental form of $\partial D$. Thus
$\mathrm{II}\succeq0$ because $D$ is convex. The boundary version of
Lemma \ref{lem:weighted-bochner} is
\begin{align}
\int_D g(L_\mu u)^2\dd\mu
&=
\int_D g
\left(
\nn{\nabla^2u}_{HS}^{2}
+\left\langle\nabla^2V\nabla u,\nabla u\right\rangle
\right)\dd\mu
\notag\\
&\quad
-\int_D
\left\langle\nabla^2g\nabla u,\nabla u\right\rangle\dd\mu
-2\int_D
\left\langle\nabla g,\nabla u\right\rangle L_\mu u\,\dd\mu
\notag\\
&\quad
+\int_{\partial D}
g\,\mathrm{II}
\left(
\nabla_{\partial D}u,\nabla_{\partial D}u
\right)e^{-V}\dd\mathcal H^{n-1}.
\label{eq:weighted-reilly}
\end{align}
Indeed, the proof of the lemma applies verbatim except when
$L_\mu|\nabla u|^2$ is integrated by parts. Tangentially
differentiating the Neumann condition gives
\[
\frac12\partial_\nu|\nabla u|^2
=
-\mathrm{II}
\left(
\nabla_{\partial D}u,\nabla_{\partial D}u
\right).
\]
This produces the last term in \eqref{eq:weighted-reilly}. The
integration of
\[
\div_\mu\left(
\left\langle\nabla g,\nabla u\right\rangle\nabla u
\right)
\]
has no boundary contribution because its normal component contains the
factor $\partial_\nu u$.

Using $L_\mu u=f$, integration by parts and
\eqref{eq:weighted-reilly}, we obtain
\begin{align*}
\int_D gf^2\dd\mu
&=
2\int_D gfL_\mu u\,\dd\mu
-\int_D g(L_\mu u)^2\dd\mu \\
&=
-2\int_D g\left\langle\nabla f,\nabla u\right\rangle\dd\mu \\
&\quad
-\int_D g
\left(
\nn{\nabla^2u}_{HS}^{2}
+\left\langle\nabla^2V\nabla u,\nabla u\right\rangle
\right)\dd\mu \\
&\quad
+\int_D
\left\langle\nabla^2g\nabla u,\nabla u\right\rangle\dd\mu \\
&\quad
-\int_{\partial D}
g\,\mathrm{II}
\left(
\nabla_{\partial D}u,\nabla_{\partial D}u
\right)e^{-V}\dd\mathcal H^{n-1}.
\end{align*}
Here the two terms containing
$f\left\langle\nabla g,\nabla u\right\rangle$ cancel. Since
$g\geq0$, $\nabla^2g\preceq0$, and $\mathrm{II}\succeq0$, it follows
that
\[
\int_D gf^2\dd\mu
\leq
\int_D g
\left(
-2\left\langle\nabla f,\nabla u\right\rangle
-\left\langle\nabla^2V\nabla u,\nabla u\right\rangle
\right)\dd\mu .
\]
For a positive definite matrix $A$ and vectors $a,b$,
\[
-2\langle a,b\rangle-\langle Ab,b\rangle
\leq
\left\langle A^{-1}a,a\right\rangle .
\]
Taking
\[
A=\nabla^2V,
\qquad
a=\nabla f,
\qquad
b=\nabla u,
\]
we conclude that
\[
\int_D gf^2\dd\mu
\leq
\int_D
g\left\langle
(\nabla^2V)^{-1}\nabla f,\nabla f
\right\rangle\dd\mu .
\]

If $\nabla^2V$ is only positive semidefinite, the same argument remains
valid provided the inverse quadratic form is understood in the extended
sense
\begin{equation}
\left\langle A^{-1}a,a\right\rangle
:=
\sup_{z\in\mathbb R^n}
\left\{
2\langle a,z\rangle-\langle Az,z\rangle
\right\}
\in[0,+\infty].
\label{eq:extended-inverse}
\end{equation}
Indeed, the preceding pointwise estimate follows by choosing $z=-b$.
When $A$ is positive definite, \eqref{eq:extended-inverse} is the usual
inverse quadratic form.

We now pass to an arbitrary convex set $K$ with
$0<\mu(K)<\infty$. Choose an increasing sequence of bounded smooth
convex domains $D_j$ such that
\[
\overline{D_j}\subset\operatorname{int}K,
\qquad
\bigcup_{j=1}^{\infty}D_j=\operatorname{int}K.
\]
Set
\[
a_j=\frac{1}{\mu(D_j)}\int_{D_j}f\,\dd\mu.
\]
Since $f\in L^1(K,\mu)$ and $\int_K f\,\dd\mu=0$, we have
$a_j\to0$.

Fix $j$. A finite concave function is continuous in the interior of
its domain. Since $\overline{D_j}\subset\operatorname{int}K$, we may
convolve $g$ with a standard mollifier on a neighborhood of
$\overline{D_j}$. In this way we obtain smooth nonnegative concave
functions $g_{j,\varepsilon}$ such that
\[
g_{j,\varepsilon}\longrightarrow g
\qquad\text{uniformly on }\overline{D_j}.
\]
If $g$ is not identically zero, its nonnegativity and concavity imply
that $g>0$ on $\operatorname{int}K$; hence $g$ has a positive minimum
on $\overline{D_j}$. Thus, whenever the right-hand side of the desired
inequality is finite, its unweighted quadratic-form integrand is
integrable on $D_j$. If the right-hand side is infinite, the conclusion
is automatic. Consequently, uniform convergence is sufficient in the
following passage to the limit. (The case $g\equiv0$ is trivial.)
Applying the smooth-domain inequality to $f-a_j$ and
$g_{j,\varepsilon}$, and then letting $\varepsilon\downarrow0$, gives
\[
\int_{D_j}g(f-a_j)^2\dd\mu
\leq
\int_{D_j}
g\left\langle
(\nabla^2V)^{-1}\nabla f,\nabla f
\right\rangle\dd\mu .
\]
Finally, Fatou's lemma on the left and monotone convergence on the
right yield
\begin{align*}
\int_K gf^2\dd\mu
&\leq
\liminf_{j\to\infty}
\int_{D_j}g(f-a_j)^2\dd\mu \\
&\leq
\lim_{j\to\infty}
\int_{D_j}
g\left\langle
(\nabla^2V)^{-1}\nabla f,\nabla f
\right\rangle\dd\mu \\
&=
\int_K
g\left\langle
(\nabla^2V)^{-1}\nabla f,\nabla f
\right\rangle\dd\mu .
\end{align*}
The boundary of a full-dimensional convex set has Lebesgue measure
zero, so replacing $\operatorname{int}K$ by $K$ does not change the
integrals. This proves the theorem.
\end{proof}

\begin{remark}
The proof in fact yields a stronger form of
Theorem~\ref{thm:BL-concave}. Assume first that $g$ is positive and of
class $C^{2}$, and set
\[
A_g=\nabla^{2}V-\frac{\nabla^{2}g}{g}.
\]
Indeed, in the weighted Bochner identity the two corresponding
quadratic terms combine as
\[
\int_K g
\left\langle \nabla^{2}V\nabla u,\nabla u\right\rangle\,d\mu
-
\int_K
\left\langle \nabla^{2}g\nabla u,\nabla u\right\rangle\,d\mu
=
\int_K g
\left\langle A_g\nabla u,\nabla u\right\rangle\,d\mu.
\]
Consequently, if $A_g$ is positive definite, then for every $f$ such
that
\[
\int_K f\,d\mu=0
\]
the same argument gives
\[
\int_K gf^{2}\,d\mu
\leq
\int_K g
\left\langle A_g^{-1}\nabla f,\nabla f\right\rangle\,d\mu.
\]
If $g$ is concave, then
\[
A_g\geq \nabla^{2}V.
\]
Thus, whenever the relevant matrices are positive definite,
\[
A_g^{-1}\leq \left(\nabla^{2}V\right)^{-1},
\]
and the preceding inequality strengthens
\eqref{eq:BL-weighted}.

More generally, concavity is not needed for the estimate involving
$A_g$: it is enough to assume that
\[
g\nabla^{2}V-\nabla^{2}g
\]
is positive definite. A corresponding degenerate statement can be
obtained by approximation when this matrix is only positive
semidefinite.
\end{remark}

\begin{remark}[Bounded domains, balls, and polytopes]
The preceding strengthening is particularly useful on bounded
domains. Let $K\subset\mathbb R^{n}$ be bounded and take $V=0$, so that
$\mu$ is normalized Lebesgue measure on $K$. Retaining the Hessian term
in the proof and using
\[
\left\|\nabla^{2}u\right\|_{\mathrm{HS}}^{2}
\geq
\frac{1}{n}(\Delta u)^{2}
\]
gives the following dimensional estimate. If $g$ is positive and
concave and $-\nabla^{2}g$ is positive definite, then
\[
\int_K gf^{2}\,dx
\leq
\frac{n}{n+1}
\int_K g^{2}
\left\langle
\left(-\nabla^{2}g\right)^{-1}\nabla f,\nabla f
\right\rangle\,dx,
\qquad
\int_K f\,dx=0.
\]
The normalization of Lebesgue measure is immaterial, since the same
constant occurs on both sides.

Equivalently, let
\[
d\nu(x)=\frac{g(x)\,dx}{\int_K g\,dx}.
\]
Applying the preceding inequality to
\[
f-\frac{1}{|K|}\int_K f\,dx
\]
and then minimizing the left-hand side over constants gives
\[
\operatorname{Var}_{\nu}(f)
\leq
\frac{n}{n+1}
\int_K
\left\langle
\left(-\frac{\nabla^{2}g}{g}\right)^{-1}
\nabla f,\nabla f
\right\rangle\,d\nu.
\]

For example, on the Euclidean unit ball, take
\[
g(x)=1-|x|^{2}.
\]
Since
\[
-\nabla^{2}g=2I,
\]
we obtain
\[
\int_{B_2^n}(1-|x|^{2})f^{2}\,dx
\leq
\frac{n}{2(n+1)}
\int_{B_2^n}(1-|x|^{2})^{2}|\nabla f|^{2}\,dx,
\qquad
\int_{B_2^n}f\,dx=0.
\]
Equivalently, if $\nu$ is the probability measure on $B_2^n$ with
density proportional to $1-|x|^{2}$, then
\[
\operatorname{Var}_{\nu}(f)
\leq
\frac{n}{2(n+1)}
\int_{B_2^n}(1-|x|^{2})|\nabla f|^{2}\,d\nu.
\]
In dimension one this inequality is sharp, with equality for
$f(x)=x$. The use of a weight vanishing on the boundary is justified
by replacing $g$ with $g+\varepsilon$ and letting
$\varepsilon\downarrow0$.

There is also a natural polytope example. Let
\[
P=
\left\{
x\in\mathbb R^{n}:
\ell_i(x)=b_i-\langle a_i,x\rangle>0,\quad
i=1,\ldots,m
\right\}
\]
be a bounded full-dimensional polytope, where the inequalities define
the facets of $P$, each facet being listed once. Set
\[
g(x)=\left(\prod_{i=1}^{m}\ell_i(x)\right)^{1/m}.
\]
The function $g$ is concave. Define
\[
v_i(x)=\frac{a_i}{\ell_i(x)},
\qquad
\overline v(x)=\frac1m\sum_{i=1}^{m}v_i(x).
\]
A direct calculation gives
\[
-\frac{\nabla^{2}g(x)}{g(x)}
=
\frac1m\sum_{i=1}^{m}v_i(x)\otimes v_i(x)
-
\overline v(x)\otimes\overline v(x).
\]
Denote this matrix by $H_P(x)$. Thus $H_P(x)$ is the covariance
matrix, with equal weights, of the vectors
\[
\frac{a_i}{\ell_i(x)}.
\]
It is a centered version of the logarithmic-barrier metric.

For every $x\in P$, the matrix $H_P(x)$ is positive definite. Indeed,
if
\[
\left\langle H_P(x)z,z\right\rangle=0,
\]
then the numbers $\langle v_i(x),z\rangle$ are all equal. If their
common value is zero, then $P$ contains the entire line $x+\mathbb Rz$.
If their common value is nonzero, all the supporting hyperplanes of
$P$ pass through a common point and $P$ contains an unbounded ray.
Both conclusions contradict the boundedness of $P$.

Consequently,
\[
\int_P gf^{2}\,dx
\leq
\frac{n}{n+1}
\int_P g
\left\langle H_P^{-1}\nabla f,\nabla f\right\rangle\,dx,
\qquad
\int_P f\,dx=0.
\]
If
\[
d\nu_P(x)=\frac{g(x)\,dx}{\int_P g\,dx},
\]
then one also has the Poincar\'e-type inequality
\[
\operatorname{Var}_{\nu_P}(f)
\leq
\frac{n}{n+1}
\int_P
\left\langle H_P^{-1}\nabla f,\nabla f\right\rangle\,d\nu_P.
\]

When $P$ is a simplex, $m=n+1$, and after an affine identification
with the standard simplex the measure $\nu_P$ is the symmetric
Dirichlet distribution whose parameters are
\[
\alpha_1=\cdots=\alpha_{n+1}
=
1+\frac{1}{n+1}
=
\frac{n+2}{n+1}.
\]
The polytope estimate follows from the smooth case by approximation
from the interior. All the preceding inequalities are affine
invariant.
\end{remark}

\section{An improved Brascamp--Lieb inequality in the unconditional case}

In this section we give the short proof of Theorem \ref{thm:poincare-uncond}.
\label{sec:unconditional}

\begin{proof}[Proof of Theorem \ref{thm:poincare-uncond}]
Write $\dd\nu=e^{-U}\dd\mu$, where $U$ is convex and unconditional. Since
$u$ is unconditional, each partial derivative $\partial_i u$ is odd in the
$i$-th variable, and we may write
\[
\partial_i u=x_i g_i
\]
with $g_i$ smooth and even in the $i$-th coordinate. Using
$T_{ij}x_ix_j\ge0$ and $ab\le(a^2+b^2)/2$,
\begin{align*}
\int\langle T\nabla u,\nabla u\rangle\dd\nu
&=
\sum_{i,j=1}^n
\int T_{ij}x_ix_j g_i g_j\dd\nu 
\le
\sum_{i,j=1}^n
\int T_{ij}x_ix_j\frac{g_i^2+g_j^2}{2}\dd\nu \\
&=
\sum_{i,j=1}^n
\int T_{ij}x_ix_j g_i^2\dd\nu 
=
\sum_{i=1}^n
\int (Tx)_i x_i g_i^2\dd\nu \\
&=
\sum_{i=1}^n
\int \partial_iV\,x_i g_i^2\dd\nu .
\end{align*}
For each $i$, integration by parts with
$\dd\nu=Z_\nu^{-1}e^{-V-U}\dd x$ gives
\[
\int \partial_iV\,x_i g_i^2\dd\nu
=
\int
\left(
g_i^2+2x_i g_i\partial_i g_i
-x_i g_i^2\partial_iU
\right)\dd\nu .
\]
Since $U$ is convex and unconditional, $x_i\partial_iU\ge0$. Hence
\[
\int \partial_iV\,x_i g_i^2\dd\nu
\le
\int
\left(
g_i^2+2x_i g_i\partial_i g_i
\right)\dd\nu
\le
\int|\nabla(x_i g_i)|^2\dd\nu .
\]
Summing over $i$ gives
\[
\int\langle T\nabla u,\nabla u\rangle\dd\nu
\le
\sum_{i=1}^n\int|\nabla\partial_i u|^2\dd\nu
=
\int\nn{\nabla^2u}_{HS}^{2}\dd\nu .
\]
\end{proof}

\begin{corollary}\label{cor:uncond-BL}
Let $\mu$ be an unconditional log-concave probability measure on $\RR^n$ with
$C^2$ density $Z^{-1}e^{-V}$. Let $T$ satisfy the assumptions of 
Theorem \ref{thm:poincare-uncond}, and assume in addition that $\nabla^2 V + T \succ 0$. 
Then for every probability measure $\nu$ with unconditional
log-concave density with respect to $\mu$, and every smooth unconditional
$f: \RR^n \to \RR$,
\[
\var_\nu f
\le
\int
\left\langle
(\nabla^2V+T)^{-1}\nabla f,\nabla f
\right\rangle\dd\nu .
\]
\end{corollary}

This corollary follows from Proposition \ref{prop:poincare-to-B}, with 
the observation  that the same proof continues to hold when the word ``even'' is replaced everywhere  with the word ``unconditional''. To avoid repetition we only prove Proposition 
\ref{prop:poincare-to-B} as stated:

\begin{proof}[Proof of Proposition \ref{prop:poincare-to-B}.]
Let $\dd\nu=\frac{1}{Z}e^{-W}\dd\mu$, where $W$ is convex and even. Let
$f$ be even and centered with respect to $\nu$. We solve
$L_\nu u=f$. Since $\nu$ and $f$ are even, $u$
may be chosen even. By the Bochner identity for the potential $V+W$,
\[
\int f^2\dd\nu
=
\int
\left(
\nn{\nabla^2u}_{HS}^{2}
+
\left\langle
(\nabla^2V+\nabla^2W)\nabla u,\nabla u
\right\rangle
\right)\dd\nu .
\]
Using \eqref{eq:genT-Poincare} and $\nabla^2W\succeq0$, we get
\[
\int f^2\dd\nu
\ge
\int
\left\langle
(\nabla^2V+T)\nabla u,\nabla u
\right\rangle\dd\nu .
\]
On the other hand,
\[
\int f^2\dd\nu=-\int\langle\nabla f,\nabla u\rangle\dd\nu, 
\]
and \eqref{eq:genT-BL} follows in the same way as in Theorem \ref{lem:BL} or Proposition \ref{prop:even-BL}.

It remains to recall the standard implication to the B-conjecture. Let
$K$ be an origin-symmetric convex body, and put $F(t)=\mu(e^tK)$. After
the change of variables $x=e^t y$, differentiation gives
\[
\frac{F'(t)}{F(t)}
=
\int
\left(n-\langle\nabla V(x),x\rangle\right)\dd\nu_t(x),
\]
where $\nu_t$ is the normalized restriction of $\mu$ to $e^tK$. A second
differentiation yields
\begin{align*}
(\log F)''(t)
&=
\var_{\nu_t}\left(\langle\nabla V,x\rangle\right) 
-\int
\left(
\langle\nabla^2Vx,x\rangle+\langle\nabla V,x\rangle
\right)\dd\nu_t .
\end{align*}
For $q(x)=\langle\nabla V(x),x\rangle$ we have
\[
\nabla q=(\nabla^2V)x+\nabla V=(\nabla^2V+T)x.
\]
Applying \eqref{eq:genT-BL} to the even function $q$ under $\nu_t$ gives
\[
\var_{\nu_t}(q)
\le
\int\left\langle(\nabla^2V+T)x,x\right\rangle\dd\nu_t
=
\int
\left(
\langle\nabla^2Vx,x\rangle+\langle\nabla V,x\rangle
\right)\dd\nu_t .
\]
Hence $(\log F)''(t)\le0$. 
\end{proof}

\section{The Brascamp--Lieb rotation trick}\label{sec:rotation}
In this section we prove Theorem \ref{thm:BL-rotation}:

\begin{proof}[Proof of Theorem \ref{thm:BL-rotation}.]
By adding a constant to $V$ we may assume that $\int e^{-V} = 1$. Fix a $C^1$ even function $f:\RR^n \to \RR$. We use the identity
\[
\var_{\mu}f=\frac{1}{2}\iint\left(f(x)-f(y)\right)^{2}e^{-V(x)-V(y)}\dd x\dd y.
\]
 Next we perform the change of variables $x=\frac{u+v}{\sqrt{2}}$, $y=\frac{u-v}{\sqrt{2}}$,
which has Jacobian $\equiv1$, and obtain
\[
\var_{\mu}f=\frac{1}{2}\iint\left(f\left(\frac{u+v}{\sqrt{2}}\right)-f\left(\frac{u-v}
{\sqrt{2}}\right)\right)^{2}e^{-V\left(\frac{u+v}{\sqrt{2}}\right)-V\left(\frac{u-v}{\sqrt{2}}\right)}
\dd u \dd v.
\]
For a fixed $u\in\RR^{n}$ we define 
\begin{align*}
f_{u}(v) &=f\left(\frac{u+v}{\sqrt{2}}\right)-f\left(\frac{u-v}{\sqrt{2}}\right) \\
V_{u}(v) &=V\left(\frac{u+v}{\sqrt{2}}\right)+V\left(\frac{u-v}{\sqrt{2}}\right).
\end{align*}
Note that $f_{u}$ is odd and $V_{u}$ is even, and so $\int f_{u}e^{-V_{u}}=0$.
Since $V_{u}$ is convex it follows from the usual Brascamp-Lieb inequality
that
\[
\var_{\mu}f=\frac{1}{2}\iint f_{u}(v)^{2}e^{-V_{u}(v)}\dd v\dd u\le\frac{1}{2}\iint\left\langle (\nabla^{2}V_{u})^{-1}\cdot\nabla f_{u},\nabla f_{u}\right\rangle e^{-V_{u}(v)}\dd v\dd u.
\]
Computing these derivatives (in $v$) we have 
\begin{align*}
\nabla f_{u}(v) & =\frac{1}{\sqrt{2}}\left(\nabla f\left(\frac{u+v}{\sqrt{2}}\right)+\nabla f\left(\frac{u-v}{\sqrt{2}}\right)\right)
\end{align*}
as well as 
\[
\nabla^{2}V_{u}(v)=\frac{1}{2}\left(\nabla^{2}V\left(\frac{u+v}{\sqrt{2}}\right)+\nabla^{2}V\left(\frac{u-v}{\sqrt{2}}\right)\right).
\]
 Therefore, performing a change of variables back to $x$ and $y$, we
obtain 
\begin{align*}
\var_{\mu}f & \le\frac{1}{2}\iint\left\langle \left(\frac{\nabla^{2}V(x)+\nabla^{2}V(y)}{2}\right)^{-1}\left(\frac{\nabla f(x)+\nabla f(y)}{\sqrt{2}}\right),\left(\frac{\nabla f(x)+\nabla f(y)}{\sqrt{2}}\right)\right\rangle \dd\mu(x)\dd\mu(y)\\
 & =\frac{1}{2}\iint\left\langle \left(\nabla^{2}V(x)+\nabla^{2}V(y)\right)^{-1}\left(\nabla f(x)+\nabla f(y)\right),\nabla f(x)+\nabla f(y)\right\rangle \dd\mu(x)\dd\mu(y).
\end{align*}
We now expand the inner product, and note that the ``cross term''
\[
\left\langle \left(\nabla^{2}V(x)+\nabla^{2}V(y)\right)^{-1}\nabla f(x),\nabla f(y)\right\rangle 
\]
 is odd in $y$ for every fixed $x$, so its integral vanishes. Using
this and the symmetry between $x$ and $y$ we obtain 
\begin{align*}
\var_{\mu}f & \le\iint\left\langle \left(\nabla^{2}V(x)+\nabla^{2}V(y)\right)^{-1}\nabla f(x),\nabla f(x)\right\rangle \dd\mu(x)\dd\mu(y)\\
 & =\int\left\langle \left(\int\left(\nabla^{2}V(x)+\nabla^{2}V(y)\right)^{-1}\dd\mu(y)\right)\cdot\nabla f(x),\nabla f(x)\right\rangle \dd\mu(x)\\
 & =\int\left\langle H_{V}\nabla f,\nabla f\right\rangle \dd\mu,
\end{align*}
as claimed. 
\end{proof}

\section{A Poincar\'e inequality for smooth norms}\label{sec:norms}

In this section we prove Theorem \ref{thm:Poincare-norms}. We begin by collecting the necessary
ingredients. The first ingredient gives us a more geometric representation for the inverse Hessian
$\left(\nabla^{2}\alpha\right)^{-1}$.
\begin{lemma} \label{lem:combine}
Let $\nn{\cdot}$ be a norm on $\RR^n$ with unit ball $K$, where $K$ is smooth and has strictly 
positive Gaussian curvature. Define $\alpha(x) = \frac{1}{2} \nn{x}^2$. Then for every $x\in\partial K$ and every $u,v\in\RR^{n}$
we have 
\[
\left\langle \left(\nabla^{2}\alpha(x)\right)^{-1}u,v\right\rangle =\left\langle x,u\right\rangle \left\langle x,v\right\rangle +\left\langle x,n_{K}\right\rangle \II^{-1}\left(P_{T_{x}\partial K}u,P_{T_{x}\partial K}v\right).
\]
 Here $n_{K}$ denotes the outer unit normal to $\partial K$ at $x$,
$P_{T_{x}\partial K}$ denotes the orthogonal projection onto the
tangent plane $T_{x}\partial K$, and $\II^{-1}$ denotes the inverse
of the second fundamental form of the hypersurface $\partial K$ at
the point $x$. 
\end{lemma}

\begin{proof}
Let us simplify our notation, writing $P$ for $P_{T_{x}\partial K}$. We will also use $S$ for the
shape operator $S:T_{x}\partial K\to T_{x}\partial K$, with the sign
convention that $Su=\partial_{u}\left(n_{K}\right)$. Finally, we
will suppress the dependence on $x$ whenever appropriate, writing
for example $\left(\nabla^{2}\alpha\right)^{-1}$ instead of $\left(\nabla^{2}\alpha(x)\right)^{-1}$.

Since $\alpha\equiv\frac{1}{2}$ on $\partial K$ we know that $\nabla\alpha=\lambda n_{K}$
for some $\lambda\in\RR$. To compute $\lambda$ we note that $\alpha$
is $2$-homogeneous and therefore 
\[
1=2\alpha=\left\langle \nabla\alpha,x\right\rangle =\lambda\left\langle n_{K},x\right\rangle ,
\]
 and so $\nabla\alpha=\frac{n_{K}}{\left\langle n_{K},x\right\rangle }$.
Again, since $\alpha$ is $2$-homogeneous, we have that 
\begin{equation}
\nabla^{2}\alpha\cdot x=\nabla\alpha=\frac{n_{K}}{\left\langle n_{K},x\right\rangle }.\label{eq:Hessian-x}
\end{equation}
 On the other hand, for all $\xi\in T_{x}\partial K$ we may use the
fact that $\partial_{\xi}\left(n_{K}\right)=S\xi$. We have 
\[
\partial_{\xi}\left(\left\langle n_{K},x\right\rangle \right)=\left\langle S\xi,x\right\rangle +\left\langle n_{K},\xi\right\rangle =\left\langle S\xi,x\right\rangle ,
\]
 and therefore 
\begin{align}
\nabla^{2}\alpha\cdot\xi & =\partial_{\xi}\left(\nabla\alpha\right)=\partial_{\xi}\left(\frac{n_{K}}{\left\langle n_{K},x\right\rangle }\right)=\frac{S\xi}{\left\langle n_{K},x\right\rangle }-\frac{\left\langle S\xi,x\right\rangle }{\left\langle n_{K},x\right\rangle ^{2}}n_{K}.\label{eq:Hessian-tangent}
\end{align}
It follows from (\ref{eq:Hessian-x}) and (\ref{eq:Hessian-tangent})
that if we define temporarily 
\[
Au=\left\langle x,u\right\rangle x+\left\langle n_{K},x\right\rangle S^{-1}Pu
\]
 then for all $u\in\RR^{n}$ we have 
\begin{align*}
\left(\nabla^{2}\alpha\right)\left(Au\right) & =\left\langle x,u\right\rangle \frac{n_{K}}{\left\langle n_{K},x\right\rangle }+\left\langle n_{K},x\right\rangle \left(\frac{Pu}{\left\langle n_{K},x\right\rangle }-\frac{\left\langle Pu,x\right\rangle }{\left\langle n_{K},x\right\rangle ^{2}}n_{K}\right)\\
 & =Pu+\left\langle u-Pu,x\right\rangle \frac{n_{K}}{\left\langle n_{K},x\right\rangle }\\
 & =Pu+\left\langle \left\langle u,n_{K}\right\rangle n_{K},x\right\rangle \frac{n_{K}}{\left\langle n_{K},x\right\rangle }=Pu+\left\langle u,n_{K}\right\rangle n_{K}=u,
\end{align*}
 which precisely means that $A=\left(\nabla^{2}\alpha\right)^{-1}$.
Therefore 
\begin{align*}
\left\langle \left(\nabla^{2}\alpha\right)^{-1}u,v\right\rangle  & =\left\langle Au,v\right\rangle =\left\langle x,u\right\rangle \left\langle x,v\right\rangle +\left\langle n_{K},x\right\rangle \left\langle S^{-1}Pu,v\right\rangle \\
 & =\left\langle x,u\right\rangle \left\langle x,v\right\rangle +\left\langle n_{K},x\right\rangle \II^{-1}\left(Pu,Pv\right)
\end{align*}
 as claimed. 
\end{proof}
Our second ingredient is a Poincaré-type inequality on $\partial K$:
\begin{lemma}
\label{lem:milman-poincare}Let $V:\RR^{n}\to\RR$ be a smooth, convex
function. Let $\mu_{\partial K}$ be the measure on $\partial K$
with density $e^{-V}\left\langle x,n_{K}\right\rangle $ with respect
to $\left.\HH^{n-1}\right|_{\partial K}$. Then for every $C^{1}$
function $f:\partial K\to\RR$ such that $\int_{\partial K}f\dd\mu_{\partial K}=0$
we have 
\[
\int_{\partial K}(n-1-\left\langle \nabla V,x\right\rangle )f^{2}\dd\mu_{\partial K}\le\int_{\partial K}\II^{-1}\left(\nabla_{\partial K}f,\nabla_{\partial K}f\right)\left\langle x,n_{K}\right\rangle \dd\mu_{\partial K}.
\]
\end{lemma}

Lemma \ref{lem:milman-poincare} is nothing more than a reformulation
of a result of Kolesnikov--Milman from \cite{KolesnikovMilman2018}. Specifically, restricting
Theorem 1.1 of \cite{KolesnikovMilman2018} to the manifold with boundary $K$ and to the case
$N=\infty$, one obtains the following:
\begin{proposition}
\label{prop:milman-original}Let $V:\RR^{n}\to\RR$ be a smooth, convex
function. Then for every $C^{1}$ function $g:\partial K\to\RR$ such
that $\int_{\partial K}ge^{-V}=0$ we have
\[
\int_{\partial K}\left(H-\left\langle \nabla V,n_{K}\right\rangle \right)g^{2}e^{-V}\le\int_{\partial K}\II^{-1}\left(\nabla_{\partial K}g,\nabla_{\partial K}g\right)e^{-V},
\]
 where $H$ denotes the mean curvature.
\end{proposition}

Here we use the convention that an integral on $\partial K$ is always
with respect to the Hausdorff measure $\left.\HH^{n-1}\right|_{\partial K}$
unless explicitly stated otherwise. Lemma \ref{lem:milman-poincare}
is obtained from Proposition \ref{prop:milman-original} using a careful
integration by parts:
\begin{proof}[Proof of Lemma \ref{lem:milman-poincare}. ]
 Given $f:\partial K\to\RR$ we define $g:\partial K\to\RR$ by $g=f\left\langle x,n_{K}\right\rangle $.
Since 
\[
\int_{\partial K}ge^{-V}=\int_{\partial K}f\left\langle x,n_{K}\right\rangle e^{-V}=\int_{\partial K}f\dd\mu_{\partial K}=0,
\]
 we deduce from Proposition \ref{prop:milman-original} that 
\begin{equation}
\int_{\partial K}\left(H-\left\langle \nabla V,n_{K}\right\rangle \right)g^{2}e^{-V}\le\int_{\partial K}\II^{-1}\left(\nabla_{\partial K}g,\nabla_{\partial K}g\right)e^{-V}.\label{eq:origin-bd-poincare}
\end{equation}

Similarly to Lemma \ref{lem:combine}, let us write $P$ for the orthogonal
projection to $T_{x}\partial K$ and $S:T_{x}\partial K\to T_{x}\partial K$
for the shape operator. Note that for $h=\left\langle x,n_{K}\right\rangle $
and $\xi\in T_{x}\partial K$ we have 
\[
\partial_{\xi}h=\left\langle \xi,n_{K}\right\rangle +\left\langle x,S\xi\right\rangle =\left\langle x,S\xi\right\rangle =\left\langle Px,S\xi\right\rangle =\left\langle SPx,\xi\right\rangle ,
\]
 so $\nabla_{\partial K}h=SPx$. Therefore 
\begin{align*}
\II^{-1}\left(\nabla_{\partial K}g,\nabla_{\partial K}g\right) & =\left\langle S^{-1}\nabla_{\partial K}\left(fh\right),\nabla_{\partial K}\left(fh\right)\right\rangle \\
 & =h^{2}\left\langle S^{-1}\nabla_{\partial K}f,\nabla_{\partial K}f\right\rangle +2hf\left\langle S^{-1}\nabla_{\partial K}h,\nabla_{\partial K}f\right\rangle +f^{2}\left\langle S^{-1}\nabla_{\partial K}h,\nabla_{\partial K}h\right\rangle \\
 & =h^{2}\II^{-1}\left(\nabla_{\partial K}f,\nabla_{\partial K}f\right)+\underbrace{h\left\langle Px,\nabla_{\partial K}\left(f^{2}\right)\right\rangle }_{\left(\star\right)}+\underbrace{f^{2}\left\langle Px,\nabla_{\partial K}h\right\rangle }_{\left(\ast\ast\right)}.
\end{align*}
 We now integrate by parts the term $\left(\star\right)$. Note that
\begin{align*}
\div_{\partial K}(Px) & =\div_{\partial K}x-\div_{\partial K}\left(hn_{K}\right)=n-1-\left\langle \nabla_{\partial K}h,n_{K}\right\rangle -h\div_{\partial K}\left(n_{K}\right)\\
 & =n-1-hH,
\end{align*}
 and therefore 
\begin{align*}
\int_{\partial K}h\left\langle Px,\nabla_{\partial K}\left(f^{2}\right)\right\rangle e^{-V} & =-\int_{\partial K}\div_{\partial K}\left(he^{-V}\cdot Px\right)f^{2}\\
 & =-\int_{\partial K}\left(\left\langle \left(\nabla_{\partial K}h-h\nabla_{\partial K}V\right),Px\right\rangle +h\left(n-1-hH\right)\right)f^{2}e^{-V}\\
 & =-\int_{\partial K}f^{2}\left\langle Px,\nabla_{\partial K}h\right\rangle e^{-V}+\int_{\partial K}f^{2}h\left\langle \nabla_{\partial K}V,Px\right\rangle e^{-V}\\
 & \qquad\qquad\qquad-\int_{\partial K}\left(n-1-hH\right)hf^{2}e^{-V}.
\end{align*}
 The first term in the expression above precisely cancels the term
$\left(\ast\ast\right)$, so the right-hand side of (\ref{eq:origin-bd-poincare})
can be written as 
\[
\int_{\partial K}\II^{-1}\left(\nabla_{\partial K}f,\nabla_{\partial K}f\right)h\dd\mu_{\partial K}+\int_{\partial K}f^{2}\left\langle \nabla_{\partial K}V,Px\right\rangle \dd\mu_{\partial K}-\int_{\partial K}\left(n-1-hH\right)f^{2}\dd\mu_{\partial K}.
\]
On the other hand, recalling that $g=hf$ the left-hand side of (\ref{eq:origin-bd-poincare})
is 
\[
\int_{\partial K}\left(H-\left\langle \nabla V,n_{K}\right\rangle \right)f^{2}h^{2}e^{-V}=\int_{\partial K}\left(H-\left\langle \nabla V,n_{K}\right\rangle \right)f^{2}h\dd\mu_{\partial K}.
\]
 The term $\int_{\partial K}Hf^{2}h\dd\mu_{\partial K}$ cancels,
and thus (\ref{eq:origin-bd-poincare}) is equivalent to 
\[
\int_{\partial K}\left[\left(n-1\right)-\left\langle \nabla_{\partial K}V,Px\right\rangle -\left\langle \nabla V,n_{K}\right\rangle h\right]f^{2}\dd\mu_{\partial K}\le\int_{\partial K}\II^{-1}\left(\nabla_{\partial K}f,\nabla_{\partial K}f\right)h\dd\mu_{\partial K}.
\]
 Finally, we observe that 
\[
\left\langle \nabla_{\partial K}V,Px\right\rangle +\left\langle \nabla V,n_{K}\right\rangle h=\left\langle P\left(\nabla V\right),Px\right\rangle +\left\langle \nabla V,n_{K}\right\rangle \left\langle x,n_{K}\right\rangle =\left\langle \nabla V,x\right\rangle ,
\]
 so the proof is complete. 
\end{proof}
The last ingredient we need is the following one-dimensional lemma,
which is precisely \cite[Lemma 8]{CERotem2023}:
\begin{lemma}
\label{lem:1d-poincare}Let $v,w:[0,\infty)\to\RR$ be continuous
functions and $C^{1}$ on $(0,\infty)$. Let $f$ be a $C^{1}$ function
on $[0,\infty)$ which is compactly supported (for simplicity) and
satisfies $f(0)=0$. Then for every $\alpha\ge0$ we have 
\[
\int_{0}^{\infty}\frac{v'}{t}f^{2}t^{\alpha}e^{-v-w}\dd t\le\int_{0}^{\infty}\left(\left(f'\right)^{2}+\alpha\cdot\left(\frac{f}{t}\right)^{2}-w'\frac{f^{2}}{t}\right)t^{\alpha}e^{-v-w}\dd t.
\]
 
\end{lemma}

With all our lemmas stated, we are ready to prove Theorem \ref{thm:Poincare-norms}.
\begin{proof}[Proof of Theorem \ref{thm:Poincare-norms}. ]
Since $\nu$ is even and more log-concave than $\mu$, we may write
\[
\dd\nu(x)=e^{-W(x)}\dd\mu(x) =e^{-v(\nn{x})-W(x)}\dd x,
\]
where $W$ is even and convex.

Using the standard approximation arguments one may assume that $f$ is compactly supported and that the
function $W$ is smooth. The co-area formula implies that for all
integrable functions $\phi$ we have 
\[
\int_{\RR^{n}}\phi(x) \dd x=\int_{\partial K}\int_{0}^{\infty}\phi(rx)r^{n-1}\dd r\dd\sigma_{K}(x)
\]
 where $\sigma_{K}$ is the measure on $\partial K$ with density $h(x)=\left\langle x,n_{K}\right\rangle $.
In particular, we have 
\[
\int\frac{v'\left(\nn{x}\right)}{\nn{x}}f^{2}\dd\nu=\int_{\partial K}\int_{0}^{\infty}\frac{v'(r)}{r}f(rx)^{2}r^{n-1}e^{-v(r)-W(rx)}\dd r\dd\sigma_{K}(x).
\]
Choosing $f_{x}(r)=f(rx)$ and $w_{x}(r)=W(rx)$, we have $f_{x}(0)=f(0)=0$
as $f$ is odd. Therefore we may use Lemma \ref{lem:1d-poincare}
and bound our integral by
\[
\int_{\partial K}\int_{0}^{\infty}\left[\underbrace{\left\langle \nabla f(rx),x\right\rangle ^{2}}_{(A)}+\underbrace{(n-1)\left(\frac{f(rx)}{r}\right)^{2}-\left\langle \nabla W(rx),x\right\rangle \frac{f^{2}(rx)}{r}}_{(B)}\right]r^{n-1}e^{-v(r)-W(rx)}\dd r\dd\sigma_{K}(x).
\]
For term $(B)$ we will exchange the order of the integrals, and use
Lemma \ref{lem:milman-poincare} for the functions $f_{r}(x)=\frac{f(rx)}{r}$
and $W_{r}(x)=W(rx)$. Since $f_{r}$ is odd and $he^{-W_{r}}$ is
even we have $\int_{\partial K}f_{r}he^{-W_{r}}=0$, and so by the
lemma 
\[
\int_{\partial K}\left(n-1-\left\langle \nabla W(rx),rx\right\rangle \right)\left(\frac{f(rx)}{r}\right)^{2}he^{-W(rx)}\le\int_{\partial K}\II^{-1}\left(\nabla_{\partial K}f_{r}(x),\nabla_{\partial K}f_{r}(x)\right)h^{2}e^{-W(rx)}.
\]
Note that $\nabla_{\partial K}f_{r}(x)=P_{x}\left(\nabla f(rx)\right)$,
where $P_{x}$ denotes the projection to $T_{x}\partial K$.
Therefore we see that our integral $\int\frac{v'\left(\nn{x}\right)}{\nn{x}}f^{2}\dd\nu$
is bounded by 
\[
\int_{\partial K}\int_{0}^{\infty}\left[\left\langle \nabla f(rx),x\right\rangle ^{2}+\II^{-1}\left(P_{x}\nabla f(rx),P_{x}\nabla f(rx)\right)h\right]r^{n-1}e^{-v(r)-W(rx)}\dd r\dd\sigma_{K}(x).
\]
 Using Lemma \ref{lem:combine} the right-hand side above is equal
to 
\[
\int_{\partial K}\int_{0}^{\infty}\left\langle \left(\nabla^{2}\alpha(rx)\right)^{-1}\nabla f(rx),\nabla f(rx)\right\rangle r^{n-1}e^{-v(r)-W(rx)}\dd r\dd\sigma_{K}(x),
\]
where we also used the fact that $\nabla^{2}\alpha$ is 0-homogeneous
and therefore $\nabla^{2}\alpha(rx)=\nabla^{2}\alpha(x)$.
Again by the co-area formula this expression is nothing more than
\[
\int_{\RR^{n}}\left\langle \left(\nabla^{2}\alpha\right)^{-1}\nabla f,\nabla f\right\rangle \dd\nu,
\]
so the proof is complete. 
\end{proof}


\begin{thebibliography}{99}

\bibitem{Artstein-AvidanetalAsymptotic2021}
S. Artstein-Avidan, A. Giannopoulos and V. D. Milman,
\emph{Asymptotic Geometric Analysis, Part I}, Mathematical Surveys and
Monographs, vol. 202, American Mathematical Society, Providence, RI, 2015.

%\bibitem{Barthe2001}
%F. Barthe, \emph{An isoperimetric result for the Gaussian measure and
%unconditional sets}, Bull. London Math. Soc. \textbf{33} (2001), 408--416.

\bibitem{BLYZ2012}
K. J. B\"or\"oczky, E. Lutwak, D. Yang and G. Zhang,
\emph{The log-Brunn--Minkowski inequality}, Adv. Math. \textbf{231}
(2012), no. 3--4, 1974--1997.

\bibitem{BobkovLedouxBrunnMinkowski2000}
S. G. Bobkov and M. Ledoux,
\emph{From Brunn--Minkowski to Brascamp--Lieb and to logarithmic Sobolev
inequalities}, Geom. Funct. Anal. \textbf{10} (2000), no. 5, 1028--1052.


\bibitem{BrascampLiebExtensions1976}
H. J. Brascamp and E. H. Lieb,
\emph{On extensions of the Brunn--Minkowski and Pr\'ekopa--Leindler
theorems, including inequalities for log concave functions, and with an
application to the diffusion equation}, J. Funct. Anal. \textbf{22}
(1976), no. 4, 366--389.

\bibitem{BrascampLieb1975}
H. J. Brascamp and E. H. Lieb,
\emph{Some inequalities for Gaussian measures and the long-range order of
the one-dimensional plasma}.
In 
\emph{Functional Integration and its Applications, Proceedings of the Conference on
Functional Integration, Cumberland Lodge, England}, edited by A. M. Arthurs, pages 1--14. 
Clarendon Press, 1975.


% \bibitem{BrascampLiebYoung1976}
% H. J. Brascamp and E. H. Lieb,
% \emph{Best constants in Young's inequality, its converse, and its
% generalization to more than three functions}, Adv. Math. \textbf{20}
% (1976), no. 2, 151--173.

\bibitem{CorderoEskenazis2025}
D. Cordero-Erausquin and A. Eskenazis, 
\emph{Concavity principles for weighted marginals},
arXiv preprint arXiv:2506.16941 (2025).


\bibitem{CEFM2004}
D. Cordero-Erausquin, M. Fradelizi and B. Maurey,
\emph{The (B) conjecture for the Gaussian measure of dilates of symmetric
convex sets and related problems}, J. Funct. Anal. \textbf{214} (2004),
no. 2, 410--427.

\bibitem{CERotem2023}
D. Cordero-Erausquin and L. Rotem,
\emph{Improved log-concavity for rotationally invariant measures of
symmetric convex sets}, Ann. Probab. \textbf{51} (2023), no. 3, 987--1003.

%\bibitem{Chernoff1981}
%H. Chernoff, \emph{A note on an inequality involving the normal
%distribution}, Ann. Probab. \textbf{9} (1981), no. 3, 533--535.

\bibitem{ENT2018}
A. Eskenazis, P. Nayar and T. Tkocz,
\emph{Gaussian mixtures: entropy and geometric inequalities},
Ann. Probab. \textbf{46} (2018), no. 5, 2908--2945.

\bibitem{Evans}
L.~C. Evans,
\emph{Partial Differential Equations},
2nd ed., Graduate Studies in Mathematics, vol.~19,
American Mathematical Society, Providence, RI, 2010.

\bibitem{GrossLogarithmic1975}
L. Gross, \emph{Logarithmic Sobolev inequalities}, Amer. J. Math.
\textbf{97} (1975), no. 4, 1061--1083.

\bibitem{KolLiv}
A.~V. Kolesnikov and G.~V. Livshyts,
\emph{On the Gardner--Zvavitch conjecture: symmetry in inequalities
of Brunn--Minkowski type},
Adv. Math. \textbf{384} (2021), Paper No.~107689, 23 pp.

\bibitem{KolesnikovMilman2017}
A. V. Kolesnikov and E. Milman,
\emph{Sharp Poincar\'e-type inequality for the Gaussian measure on the
boundary of convex sets}, in \emph{Geometric Aspects of Functional
Analysis}, Lecture Notes in Math., vol. 2169, Springer, Cham, 2017,
pp. 221--234.

\bibitem{KolesnikovMilman2018}
A. V. Kolesnikov and E. Milman,
\emph{Poincar\'e and Brunn--Minkowski inequalities on the boundary of weighted Riemannian manifolds
}, 
Amer. J. Math. \textbf{140} (2018), no. 5, 1147--1185. 

\bibitem{Latala2002}
R. Lata\l{}a, \emph{On some inequalities for Gaussian measures},
Proceedings of the International Congress of Mathematicians, Vol. II
(Beijing, 2002), Higher Ed. Press, Beijing, 2002, pp. 813--822.

\bibitem{LeindlerCertain1972}
L. Leindler, \emph{On a certain converse of H\"older's inequality. II},
Acta Sci. Math. (Szeged) \textbf{33} (1972), 217--223.

\bibitem{openAIlogBM2026}
OpenAI, \emph{The logarithmic Brunn--Minkowski conjecture},
OpenAI Math Release preprint
\url{https://github.com/openai/math/blob/main/preprints/The-logarithmic-Brunn-Minkowski-conjecture-September-23-2026/paper.pdf} (2026)


\bibitem{PrekopaLogarithmic1971}
A. Pr\'ekopa, \emph{Logarithmic concave measures with application to
stochastic programming}, Acta Sci. Math. (Szeged) \textbf{32} (1971),
301--316.

\bibitem{PrekopaLogarithmic1973}
A. Pr\'ekopa, \emph{On logarithmic concave measures and functions},
Acta Sci. Math. (Szeged) \textbf{34} (1973), 335--343.

\bibitem{Saroglou2016}
C. Saroglou, \emph{More on logarithmic sums of convex bodies},
Mathematika \textbf{62} (2016), no. 3, 818--841.

\end{thebibliography}
\end{document}